\documentclass[12pt,amssymb,amstex]{article}

\usepackage{amsmath}
\usepackage{amssymb}
\usepackage{amsxtra}
\usepackage{latexsym}
\usepackage{ifthen}

\usepackage{mathtext}
\usepackage[cp1251]{inputenc}
\usepackage[T2A]{fontenc}
\usepackage[dvips]{graphicx}

\begin{document}

\newcounter{lemma}[section]
\newcommand{\lemma}{\par \refstepcounter{lemma}%
{\bf Lemma \arabic{section}.\arabic{lemma}.}}
\renewcommand{\thelemma}{\thesection.\arabic{lemma}}

\newcounter{corol}[section]
\newcommand{\corol}{\par \refstepcounter{corol}%
{\bf Corollary \arabic{section}.\arabic{corol}.}}
\renewcommand{\thecorol}{\thesection.\arabic{corol}}

\newcounter{rem}[section]
\newcommand{\rem}{\par \refstepcounter{rem}%
{\bf Remark \arabic{section}.\arabic{rem}.}}
\renewcommand{\therem}{\thesection.\arabic{rem}}

\newcounter{theo}[section]
\newcommand{\theo}{\par \refstepcounter{theo}%
{\bf Theorem \arabic{section}.\arabic{theo}.}}
\renewcommand{\thetheo}{\thesection.\arabic{theo}}

\newcounter{propo}[section]
\newcommand{\propo}{\par \refstepcounter{propo}%
{\bf Proposition \arabic{section}.\arabic{propo}.}}
\renewcommand{\thepropo}{\thesection.\arabic{propo}}

\numberwithin{equation}{section}

\newcommand{\osc}{\operatornamewithlimits{osc}}

\def\Xint#1{\mathchoice
   {\XXint\displaystyle\textstyle{#1}}%
   {\XXint\textstyle\scriptstyle{#1}}%
   {\XXint\scriptstyle\scriptscriptstyle{#1}}%
   {\XXint\scriptscriptstyle\scriptscriptstyle{#1}}%
   \!\int}
\def\XXint#1#2#3{{\setbox0=\hbox{$#1{#2#3}{\int}$}
     \vcenter{\hbox{$#2#3$}}\kern-.5\wd0}}
\def\dashint{\Xint-}

\markboth{\centerline{D. KOVTONYUK, V. RYAZANOV, R. SALIMOV, E.
SEVOST'YANOV}} {\centerline{ON MAPPINGS IN THE ORLICZ-SOBOLEV
CLASSES}}

\def\cc{\setcounter{equation}{0}
\setcounter{figure}{0}\setcounter{table}{0}}

\overfullrule=0pt

\title{{\bf On mappings in the Orlicz-Sobolev classes}}

\author{{\bf D. Kovtonyuk, V. Ryazanov, R. Salimov, E. Sevost'yanov}\\
{}\\ {IN MEMORY OF ALBERTO CALDERON (1920--1998)}}

\date{\today \hskip 4mm ({\tt OS-120111-ARXIV.tex})}
\maketitle

\Large \abstract  First of all, we prove that open mappings in
Orlicz-Sobolev classes $W^{1,\varphi}_{\rm loc}$ under the Calderon
type condition on $\varphi$ have the total differential a.e. that is
a generalization of the well-known theorems of
Gehring-Lehto-Menchoff in the plane and of V\"ais\"al\"a in ${\Bbb
R}^n$, $n\geqslant3$. Under the same condition on $\varphi$, we show
that continuous mappings $f$ in $W^{1,\varphi}_{\rm loc}$, in
particular, $f\in W^{1,p}_{\rm loc}$ for $p>n-1$ have the
$(N)$-property by Lusin on a.e. hyperplane. Our examples demonstrate
that the Calderon type condition is not only sufficient but also
necessary for this and, in particular, there exist homeomorphisms in
$W^{1,n-1}_{\rm loc}$ which have not the $(N)$-property with respect
to the $(n-1)$-dimensional Hausdorff measure on a.e. hyperplane. It
is proved on this base that under this condition on $\varphi$ the
homeomorphisms $f$ with finite distortion in $W^{1,\varphi}_{\rm
loc}$ and, in particular, $f\in W^{1,p}_{\rm loc}$ for $p>n-1$ are
the so-called lower $Q$-homeomorphisms where $Q(x)$ is equal to its
outer dilatation $K_f(x)$ as well as the so-called ring
$Q_*$-homeomorphisms with $Q_*(x)=\left[K_{f}(x)\right]^{n-1}$. This
makes possible to apply our theory of the local and boundary
behavior of the lower and ring $Q$-homeomorphisms to homeomorphisms
with finite distortion in the Orlicz-Sobolev classes. \endabstract

\bigskip



\bigskip

\begin{center}{\Large\bf Contents}\end{center}

\begin{itemize}

\item[1.] Introduction \dotfill \pageref{1}

\item[2.] Preliminaries \dotfill \pageref{2}

\item[3.] Differentiability of open mappings \dotfill \pageref{3}

\item[4.] The Lusin and Sard properties on surfaces \dotfill \pageref{4}

\item[5.] On BMO and FMO functions \dotfill \pageref{5}

\item[6.] On some integral conditions \dotfill \pageref{6}

\item[5.] Moduli of families of surfaces \dotfill \pageref{7}

\item[6.] Lower and ring $Q$-homeomorphisms \dotfill \pageref{8}

\item[7.] Lower $Q$-homeomorphisms and Orlicz-Sobolev classes
\dotfill \pageref{9}

\item[8.] Equicontinuous and normal families \dotfill \pageref{10}

\item[9.] On domains with regular boundaries \dotfill \pageref{11}

\item[12.] The boundary behavior \dotfill \pageref{12}

\item[13.] Some examples \dotfill \pageref{13}

\end{itemize}

\Large

\section{Introduction}\label{1}

Moduli provide us the main geometric tool in the mapping theory. The
recent development of the moduli method in the connection with
modern classes of mappings can be found in the monograph
\cite{MRSY}, see also recent books in the moduli and capacity theory
\cite{BBZ}, \cite{Dub} and
 \cite{Vas} as well as the following papers and
monographs \cite{AIM}, \cite{Koval}, \cite{Kuz}, \cite{Shl},
\cite{Sol}, \cite{Sych}, \cite{Zo3} and further references therein.
In the present paper we show that the theories of the so-called
lower and ring $Q$-homeomorphisms developed in \cite{MRSY} can be
applied to a wide range of mappings with finite distortion in the
Orlicz-Sobolev classes. Note that the plane case has been recently
studied in \cite{KPR} and \cite{LSS}. Recall, it was established
therein that each homeomorphism of finite distortion in the plane is
a lower and ring $Q$-homeomorphism with $Q(x)=K_f(x)$.

In what follows, $D$ is a domain in a finite-dimensional Euclidean
space. Following Orlicz, see \cite{Or1}, given a convex increasing
function $\varphi:[0,\infty)$ $\to[0,\infty)$, $\varphi(0)=0$,
denote by $L_{\varphi}$ the space of all functions $f:D\to{\Bbb R}$
such that
\begin{equation}\label{eqOS1.1}\int\limits_{D}\varphi\left(\frac{|f(x)|}
{\lambda}\right)\,dm(x)<\infty\end{equation} for some $\lambda>0$
where $dm(x)$ corresponds to the Lebesgue measure in $D$, see also
the monographs \cite{KR} and \cite{Za}. $L_{\varphi}$ is called the
{\bf Orlicz space}.  If $\varphi(t)=t^p$, then we write $L_{p}$. In
other words, $L_{\varphi}$ is the cone over the class of all
functions $g:D\to{\Bbb R}$ such that
\begin{equation}\label{eqOS1.2}
\int\limits_{D}\varphi\left(|g(x)|\right)\,dm(x)<\infty\end{equation}
which is also called the {\bf Orlicz class}, see \cite{BO}.
\medskip

The {\bf Orlicz-Sobolev class} $W^{1,\varphi}_{\rm loc}(D)$ is the
class of all locally integrable functions $f$ given in $D$ with the
first distributional de\-ri\-va\-ti\-ves whose gradient $\nabla f$
belongs locally in $D$ to the Orlicz class. Note that by definition
$W^{1,\varphi}_{\rm loc}\subseteq W^{1,1}_{\rm loc}$. As usual, we
write $f\in W^{1,p}_{\rm loc}$ if $\varphi(t)=t^p$, $p\geqslant1$.
It is known that a continuous function $f$ belongs to $W^{1,p}_{\rm
loc}$ if and only if $f\in ACL^{p}$, i.e., if $f$ is locally
absolutely continuous on a.e. straight line which is parallel to a
coordinate axis and if the first partial derivatives of $f$ are
locally integrable with the power $p$, see, e.g., 1.1.3 in
\cite{Maz}. The concept of the distributional derivative was
introduced by Sobolev \cite{So} in ${\Bbb R}^n$, $n\geqslant2$, and
now it is developed under wider settings, see, e.g., \cite{Am},
\cite{Ha}, \cite{He}, \cite{HK$^*_2$}, \cite{HKST}, \cite{Ma$_1$},
\cite{MRSY}, \cite{Re$_1$}, \cite{UV} and \cite{UkhVo}.
\medskip

Later on, we also write $f\in W^{1,\varphi}_{\rm loc}$ for a locally
integrable vector-function $f=(f_1,\ldots,f_m)$ of $n$ real
variables $x_1,\ldots,x_n$ if $f_i\in W^{1,1}_{\rm loc}$ and
\begin{equation}\label{eqOS1.2a}
\int\limits_{D}\varphi\left(|\nabla
f(x)|\right)\,dm(x)<\infty\end{equation} where $|\nabla
f(x)|=\sqrt{\sum\limits_{i,j}\left(\frac{\partial f_i}{\partial
x_j}\right)^2}$. In the main part of the paper we use the notation
$W^{1,\varphi}_{\rm loc}$ for more general functions $\varphi$ than
in the classical Orlicz classes giving up the condition on convexity
of $\varphi$. In fact we need the convexity of $\varphi$ only in
Section 13. Note that the Orlicz-Sobolev classes are intensively
studied in various aspects, see, e.g., \cite{AC}, \cite{AIM},
\cite{Bi$_2$}, \cite{Ca}, \cite{Ci}, \cite{Do}, \cite{GR},
\cite{Hs}, \cite{IKO}, \cite{KP}, \cite{Ko}, \cite{LM}, \cite{LL},
\cite{Tu} and \cite{Vu}.
\medskip

Recall that a homeomorphism $f$ between domains $D$ and $D'$ in
${\Bbb R}^n$, $n\geqslant2$, is called of {\bf finite distortion}
if $f\in W^{1,1}_{\rm loc}$ and
\begin{equation}\label{eqOS1.3}\Vert
f'(x)\Vert^n\leqslant K(x)\cdot J_f(x)\end{equation} with a.e.
finite function $K$ where $\Vert f'(x)\Vert$ denotes the matrix norm
of the Jacobian matrix $f'$ of $f$ at $x\in D$,
$||f'(x)||=\sup\limits_{h\in{\Bbb R}^n,|h|=1}|f'(x)\cdot h|$, and
$J_f(x)={\rm det}f'(x)$ is its Jacobian. Later on, we use the
notation $K_{f}(x)$ for the minimal function $K(x)\geqslant1$ in
(\ref{eqOS1.3}), i.e., we set $K_f(x)=\Vert f'(x)\Vert^n/$ $J_f(x)$
if $J_f(x)\ne 0,$ $K_f(x)=1$ if $f'(x)=0$ and $K_f(x)=\infty$ at the
rest points.
\medskip

First this notion was introduced on the plane for $f\in W^{1,2}_{\rm
loc}$ in the work \cite{IS}. Later on, this condition was changed by
$f\in W^{1,1}_{\rm loc}$ but with the additional condition $J_f\in
L^1_{\rm loc}$ in the monograph \cite{IM}. The theory of the
mappings with finite distortion had many successors, see, e.g.,
\cite{AIKM}, \cite{Cr1}, \cite{Cr2}, \cite{FKZ}, \cite{GI},
\cite{HK$^*_1$}--\cite{HM$^*$}, \cite{HP}, \cite{IKO$_1$},
\cite{IKO}, \cite{Ka}, \cite{KKM$_1$}, \cite{KKM$_2$}, \cite{KKMOZ},
\cite{KO$_1$}--\cite{KR$_*$}, \cite{MV$_1$}, \cite{MV$_2$},
\cite{On$_1$}--\cite{On$_3$}, \cite{Pa} and
\cite{Ra$_1$}--\cite{Ra$_4$}. They had as predecessors the mappings
with bounded distortion, see \cite{Re} and \cite{Vo}, in other
words, the quasiregular mappings, see, e.g., \cite{HKM}, \cite{MRV}
and \cite{Ri}. They are also closely related to the earlier mappings
with the bounded Dirichlet integral, see, e.g., \cite{LF} and
\cite{Su1}--\cite{Su3}, and the mappings quasiconformal in the mean
which had a rich history, see, e.g., \cite{Ah}, \cite{Bel},
\cite{Bi$_1$}, \cite{Go}, \cite{GK}, \cite{GMSV}, \cite{GRSY},
\cite{Kr1}--\cite{Ku2}, \cite{Per1}, \cite{Per2}, \cite{Pes},
\cite{Rya1}--\cite{Rya3}, \cite{RSY$_2$}, \cite{Str}, \cite{StrSy},
\cite{UkhVo}, \cite{Zo1} and \cite{Zo2}.
\medskip

Note that the above additional condition $J_f\in L^1_{\rm loc}$ in
the definition of the mappings with finite distortion can be omitted
for ho\-meo\-mor\-phisms. Indeed, for each homeomorphism $f$ between
domains $D$ and $D'$ in ${\Bbb R}^n$ with the first partial
derivatives a.e. in $D$, there is a set $E$ of the Lebesgue measure
zero such that $f$ satisfies $(N)$-property by Lusin on $D\setminus
E$ and
\begin{equation}\label{eqOS1.1.1}\int\limits_{A}J_f(x)\,dm(x)=|f(A)|\end{equation}
for every Borel set $A\subset D\setminus E$, see, e.g., 3.1.4, 3.1.8
and 3.2.5 in \cite{Fe}. On this base, it is easy by the H\"older
inequality to verify, in particular, that if $f\in W^{1,1}_{\rm
loc}$ is a ho\-meo\-mor\-phism and $K_f\in L^q_{\rm loc}$ for
$q>n-1$, then also $f\in W^{1,p}_{\rm loc}$ for $p>n-1$, that we
often use further to obtain corollaries.

\newpage

In this paper $H^k$, $k\geqslant 0$, denotes the {\bf $\bf
k$-dimensional Hausdorff measure} in ${\Bbb R}^n$, $n\geqslant1$.
More precisely, if $A$ is a set in ${\Bbb R}^n$, then
\begin{eqnarray}\label{eq8.2.1}H^k(A)\ &=&\ \sup_{\varepsilon>0}\
H^k_{\varepsilon}(A)\,, \\ \label{eq8.2.2} H^k_{\varepsilon}(A)\
&=&\ \inf\ \sum^{\infty}_{i=1}\left({\rm diam}\,
A_i\right)^k\,,\end{eqnarray} where the infimum in (\ref{eq8.2.2})
is taken over all coverings of $A$ by sets $A_i$ with ${\rm
diam}\,A_i<\varepsilon$, see, e.g., \cite{Mattila} in this
connection. Note that $H^k$ is an {\bf outer measure in the sense of
Caratheodory}, i.e.,
\medskip

\noindent (1) $\ H^k(X)\ \leqslant\ H^k(Y)\ $ whenever $X\subseteq
Y$,
\medskip

\noindent (2) $\ H^k(\bigcup\limits_i X_i)\ \leqslant\
\sum\limits_i H^k(X_i)\ $ for each sequence of sets $X_i$,
\medskip

\noindent (3) $\ H^k(X\cup Y)\ =\ H^k(X)+H^k(Y)\ $ whenever $\mathrm
{dist}(X,Y)>0$.
\bigskip

A set $E\subset{\Bbb R}^n$ is called {\bf measurable} with respect
to $H^k$ if $H^k(X)=H^k(X\cap E)+H^k(X\setminus E)$ for every set
$X\subset{\Bbb R}^n$. It is well known that every Borel set is
measurable with respect to any outer measure in the sense of
Caratheodory, see, e.g., Theorem II (7.4) in \cite{Sa}. Moreover,
$H^k$ is Borel regular, i.e., for every set $X\subset{\Bbb R}^n$,
there is a Borel set $B\subset{\Bbb R}^n$ such that $X\subset B$ and
$H^k(X)=H^k(B)$, see, e.g., Theorem II (8.1) in \cite{Sa} and
Section 2.10.1 in \cite{Fe}. The latter implies that, for every
measurable set $E\subset{\Bbb R}^n$, there exist Borel sets $B_*$
and $B^*\subset{\Bbb R}^n$ such that $B_*\subset E\subset B^*$ and
$H^k(B^*\setminus B_*)=0$, see, e.g., Section 2.2.3 in \cite{Fe}. In
particular, $H^k(B^*)=H^k(E)=H^k(B_*)$. \medskip

If $H^{k_1}(A)<\infty$, then $H^{k_2}(A)=0$ for every $k_2>k_1$,
see, e.g., VII.1.B in \cite{HW}. The quantity $$\mbox{dim}_H\, A =
\sup\limits_{H^k(A)>0}k$$ is called the {\bf Hausdorff dimension} of
a set $A$.
\bigskip

It is known that the outer Lebesgue measure
$m(A)=\Omega_n\cdot2^{-n}H^n(A)$ for sets $A$ in ${\Bbb R}^n$ where
$\Omega_n$ denotes the volume of the unit ball in ${\Bbb R}^n$, see
\cite{Sard}.
\bigskip

It was shown in \cite{GV} that a set $A$ with $\mbox{dim}_H\, A =p$
can be transformed into a set $B=f(A)$ with $\mbox{dim}_H\, B =q$
for each pair of numbers $p$ and $q\in (0,n)$ under a quasiconformal
mapping $f$ of ${\Bbb R}^n$ onto itself, cf. also \cite{Ba} and
\cite{Bi$_2$}.

\cc
\section{Preliminaries}\label{2}

First of all, the following fine property of functions $f$ in the
Sobolev classes $W^{1,p}_{\rm loc}$ was proved in the monograph
\cite{GR}, Theorem 5.5, and can be extended to the Orlicz-Sobolev
classes. The statement follows directly from the Fubini theorem and
the known characterization of functions in Sobolev's class
$W^{1,1}_{\rm loc}$ in terms of ACL (absolute continuity on lines),
see, e.g., Section 1.1.3 in \cite{Maz}, and comments in
In\-tro\-duc\-tion.

\bigskip

\begin{propo}\label{prOS2.1} {\it Let $U$ be an open set in ${\Bbb R}^n$ and
let $f:U\to{\Bbb R}^m,$ $m=1,2,\ldots$, be a mapping in the
Orlicz-Sobolev class $W^{1,\varphi}_{\rm loc}(U)$ with an
increasing function $\varphi:[0,\infty)\to[0,\infty)$. Then, for
every $k$-dimensional direction $\Gamma$ for a.e. $k$-dimensional
plane ${\mathcal P}\in\Gamma$, $k=1,2,\ldots,n-1$, the restriction
of the function $f$ on the set ${\mathcal P}\cap U$ is a function
in the class $W^{1,\varphi}_{\rm loc}({\mathcal P}\cap U)$.}
\end{propo}

\bigskip

Here the class $W^{1,\varphi}_{\rm loc}$ is well defined on a.e.
$k$-dimensional plane because partial derivatives are Borel
functions and, moreover, Sobo\-lev's classes are invariant with
respect to quasi-isometric transforma\-tions of systems of
coordinates, in particular, with respect to rotations, see, e.g.,
1.1.7 in \cite{Maz}. Recall also that a $k$-dimensional direction
$\Gamma$ in ${\Bbb R}^n$ is the class of equivalence of all
$k$-dimensional planes in ${\Bbb R}^n$ that can be obtained each
from other by a parallel shift. Note that each $(n-k)$-dimensional
plane ${\mathcal T}$ which is quite orthogonal to a $k$-dimensional
plane ${\mathcal P}$ in $\Gamma$ intersects ${\mathcal P}$ at a
single point $X({\mathcal P})$. If $E$ is a subset of $\Gamma$, then
$X(E)$ denotes the collection of all point $X({\mathcal P})$,
${\mathcal P}\in E$. It is clear that $(n-k)$-dimensional measure of
the set $X(E)$ does not depend of the choice of the plane ${\mathcal
T}$ and it is denoted by $\mu_{n-k}(E)$. They say that a property
holds for almost every plane in $\Gamma$ if $\mu_{n-k}(E)=0$ for a
set $E$ of all planes ${\mathcal P}$ for which the property fails.

\bigskip

Recall also the little--known Fadell theorem in \cite{Fa} that makes
possible us to extend the well-known theorems of
Gehring-Lehto-Menchoff in the plane and V\"ais\"al\"a in ${\Bbb
R}^n$, $n\geqslant3$, see, e.g., \cite{GL}, \cite{Me} and
\cite{Va$_2$}, on differentiability a.e. of open mappings in
Sobolev's classes to the open mappings in Orlicz-Sobolev classes in
${\Bbb R}^{n}$, $n\geqslant3$.

\bigskip

\begin{propo}\label{prOS2.0} {\it Let $f:\Omega\to{\Bbb R}^{n}$ be a
continuous open mapping on an open set $\Omega$ in $\Bbb{R}^n$,
$n\geqslant3$. If $f$ has a total differential a.e. on $\Omega$
with respect to $n-1$ variables, then $f$ has a total differential
a.e. on $\Omega$.}
\end{propo}

\bigskip

Now, let us give the following Calderon result in \cite{Ca}, p.
208, cf. Lemma 3.2 in \cite{KKM}.

\bigskip

\begin{propo}\label{prOS2.4} {\it Let $\varphi:[0,\infty)\to[0,\infty)$
be an increasing functi\-on such that $\varphi(0)=0$ and
\begin{equation}\label{eqOS2.5} A\ \colon=\ \int\limits_{0}^{\infty}\left[\frac{t}{\varphi(t)}\right]^
{\frac{1}{k-1}}\ dt\ <\ \infty\end{equation} for a natural number
$k\geqslant2$ and let $f:G\to{\Bbb R}$ be a function given in a
domain $G\subset{\Bbb R}^k$ of the class $W^{1,\varphi}(G)$. Then
\begin{equation}\label{eqOS2.6}{\rm diam}\,
f(C)\leqslant\alpha_{k}A^{\frac{k-1}{k}}\left[\
\int\limits_{C}\varphi\left(|\nabla
f|\right)\,dm(x)\right]^{\frac{1}{k}}\end{equation} for every cube
$C\subset G$ whose adges are oriented along coordinate axes where
$\alpha_k$ is a constant depending only on $k$.} \end{propo}

\bigskip

Perhaps, the Calderon work \cite{Ca} had time to be forgotten
because it was published long ago in a badly accessible journal.

\bigskip

\begin{rem}\label{remOS2.1}
It is clear that the behavior of the function $\varphi$ about $0$
is not essential and (\ref{eqOS2.6}) holds with the replacement of
the constant $A$ by the constant \begin{equation}\label{eqOS2.55}
A_*\ \colon=\ \left[{\frac{1}{\varphi(1)}}\right]^{\frac{1}{k-1}}\
+\ \int\limits_{1}^{\infty}
\left[\frac{t}{\varphi(t)}\right]^{\frac{1}{k-1}}\,dt\ <\ \infty
\end{equation} and $\varphi(t)$ by $\varphi_*(t)\equiv \varphi(1)$
for $t\in(0,1)$, $\varphi_*(0)=0$ and $\varphi_*(t)=\varphi(t)$ for
$t\geqslant1$. Indeed, applying Proposition \ref{prOS2.4} to the one
parameter family of the functions
$\varphi_{\lambda}(t)=\varphi(t)+\lambda\cdot[\varphi_*(t)-\varphi(t)]$,
$\lambda\in[0,1)$, we obtain (\ref{eqOS2.6}) with the changes
$A\mapsto A_*$ and $\varphi\mapsto \varphi_*$ as $\lambda\to1$.
\end{rem}

\bigskip

Finally, one more statement of Calderon in \cite{Ca}, p. 209,
211-212, will be also useful later on.

\bigskip

\begin{propo}\label{prOS2.5} {\it Let $\varphi:[0,\infty)\to[0,\infty)$
be a convex increasing function such that $\varphi(0)=0$ and
\begin{equation}\label{eqOS2.66} \int\limits_{1}^{\infty}\left[\frac{t}{\varphi(t)}\right]^
{\frac{1}{k-1}}dt=\infty\end{equation} for a natural number
$k\geqslant2$. Then there is a continuously differentiab-le
decreasing function $F:(0,\infty)\to[0,\infty)$ with a compact
support such that $F(t)\to\infty$ as $t\to0$, $F'(t)$ is
non-decreasing, $F'(t)\to-\infty$ as $t\to0$ and $F_*(x)=F(|x|)$,
$x\in{\Bbb R}^k$, belongs to the class $W^{1,\varphi}({\Bbb R}^k)$,
i.e., $f\in W^{1,1}({\Bbb R}^k)$ and
\begin{equation}\label{eqOS2.7} \int\limits_{{\Bbb
R}^k}\varphi\left(|\nabla F_*|\right)\,dm(x)\leqslant1.
\end{equation}} \end{propo}

\bigskip

\begin{rem}\label{remOS2.2}
The function $F$ from Proposition \ref{prOS2.5} can be describ\-ed
in a more constructive way. More precisely, set
\begin{equation}\label{eqOS2.555} \Phi(t)\ =\
\int\limits_{1}^{t}
\left[\frac{\tau}{\varphi(\tau)}\right]^{\frac{1}{k-1}}\ d\tau
\end{equation} and \begin{equation}\label{eqOS2.33} \Psi(t)\ =\
\frac{\Phi '(t)}{\Phi(t)}\ =\
\left[\frac{t}{\varphi(t)}\right]^{\frac{1}{k-1}}\frac{1}{\Phi(t)}.
\end{equation} The function $\Psi$ is continuous and decreasing and tends
to $0$ as $t\to\infty$ and to $\infty$ as $t\to1$. thus, its
inverse function $h(s)$ is well defined for all $s>0$. It was
proved by Calderon in \cite{Ca} that
\begin{equation}\label{eqOS2.333} \int\limits_{0}^{1}h(s)\,ds\ =\
\infty\,,\qquad\int\limits_{0}^{1}\varphi(h(s))\,s^{k-1}\,ds<\infty\,.
\end{equation} Then we may set $F(t)\equiv0$ for $t\geqslant1$ and
\begin{equation}\label{eqOS2.44} F(t)=\int\limits_{t}^{1}[h(s)-h(1)]\,ds\qquad
\forall\ t\in [0,1]\,.\end{equation}

On the base of Proposition \ref{prOS2.4}, it was proved by Calderon
that each continuous function $f:G\to{\Bbb R}$ in $W^{1,\varphi}(G)$
under the condition \begin{equation}\label{eqOS2.55}
\int\limits_{1}^{\infty}\left[\frac{t}{\varphi(t)}\right]^
{\frac{1}{k-1}}dt<\infty\end{equation} has a total differential a.e.
Moreover, on the base of Proposition \ref{prOS2.5}, under the
condition (\ref{eqOS2.66}) Calderon has constructed a continuous
function $f:{\Bbb R}^k\to{\Bbb R}$ which has not a total
differential a.e. We use Propositions \ref{prOS2.4} and
\ref{prOS2.5} for other purposes.
\end{rem}

\cc
\section{The differentiability of open mappings}\label{3}

Let us start from the following statement which is due to Calderon
\cite{Ca} but we prefer in comparison with \cite{Ca} to prove it on
the base of the Stepanoff theorem.

\bigskip

\begin{lemma}{}\label{lemOS3.0} {\it Let $\Omega$ be an open set  in ${\Bbb R}^k$,
$k\geqslant2$, and let $f:\Omega\to{\Bbb R}^{m}$, $m\ge 1$, be a
continuous mapping in the class $W^{1,\varphi}_{\rm loc}(\Omega)$
with some increasing $\varphi:[0,\infty)\to[0,\infty)$ such that
$\varphi(0)=0$ and
\begin{equation}\label{eqOS3.001}A:=\int\limits_{1}^{\infty}\left[\frac{t}{\varphi(t)}\right]^
{\frac{1}{k-1}}dt<\infty.\end{equation} Then $f$ has a total
differential a.e. in $\Omega$.} \end{lemma}

\bigskip

{\it Proof.} Given $x\in\Omega$, we set
$$L(x,f)=\limsup\limits_{y\to x}\frac{|f(y)-f(x)|}{|y-x|}\,.$$
By the Stepanoff theorem, see, e.g., Theorem 3.1.9 in \cite{Fe}, the
proof is reduced to the proof of the fact $L(x,f)<\infty$ a.e. in
$\Omega$.\bigskip

Denote by $C(x,r)$ the oriented cube centered at $x$ such that the
ball $B(x,r)$ is inscribed in $C(x,r)$ where $r=|x-y|$. Then
$$L(x,f)\ =\ \limsup\limits_{y\to x}\,\frac{|f(x)-f(y)|}{|x-y|}\ \leqslant$$
$$\leqslant\limsup\limits_{r\to0}\,\frac{d(fB(x,r))}{r}\
\leqslant\ \limsup\limits_{r\to0}\,\frac{d(fC(x,r))}{r}$$ and by
Proposition \ref{prOS2.4} and Remark \ref{remOS2.1} we get
$$L(x,f)\leqslant\gamma_{k,m}\,A^{\frac{k-1}{k}}_*\limsup\limits_{r\to0}
\left[\frac{1}{r^k}\int\limits_{C(x,r)}\varphi_*\left(|\nabla
f|\right)\,dm(x)\right]^{\frac{1}{k}}<\infty$$ for a.e. $x\in\Omega$
by the Lebesgue theorem on differentiability of indefini\-te
integral, see, e.g., Theorem IV.5.4 in \cite{Sa}. The proof is
complete.

\bigskip

Combining Lemma \ref{lemOS3.0} and Proposition \ref{prOS2.1}, we
obtain the following statement.

\bigskip

\begin{corol}{}\label{corOS4.2a} {\it Let $\Omega$ be an open set in ${\Bbb R}^{n}$,
$n\geqslant3$, and let $f:\Omega\to{\Bbb R}^{m}$, $m\ge 1$, be a
continuous mapping in the class $W^{1,\varphi}_{\rm loc}(\Omega)$
with an increasing function $\varphi:[0,\infty)\to[0,\infty)$ such
that $\varphi(0)=0$ and
\begin{equation}\label{eqOS3.0a}\int\limits_{1}^{\infty}\left[\frac{t}{\varphi(t)}\right]^
{\frac{1}{n-2}}dt<\infty.\end{equation} Then $f:\Omega\to{\Bbb R}^m$
has a total differential a.e. on a.e. hyperplane which is parallel
to a coordinate hyperplane.} \end{corol}

\bigskip

Combing Corollary \ref{corOS4.2a} and the Fadell result in
\cite{Fa}, see Proposition \ref{prOS2.0} above, we obtain the main
result of this section.

\bigskip

\begin{theo}{}\label{thOS4.0} {\it Let $\Omega$ be an open set in ${\Bbb R}^{n}$,
$n\geqslant3$, and let $f:\Omega\to{\Bbb R}^{n}$ be a continuous
open mapping in the class $W^{1,\varphi}_{\rm loc}(\Omega)$ with an
increasing $\varphi:[0,\infty)\to[0,\infty)$ such that
$\varphi(0)=0$ and (\ref{eqOS3.0a}) holds. Then $f$ has a total
differential a.e. in $\Omega$.} \end{theo}

\bigskip

\begin{corol}{}\label{corOS4.2b} {\it If $f:\Omega\to{\Bbb R}^{n}$ is a homeomorphism in
$W^{1,1}_{\rm loc}$ with $K_f\in L^p_{\rm loc}$ for $p>n-1$, then
$f$ is differentiable a.e. }\end{corol}

\bigskip

\begin{rem}\label{remOS4.a*} In particular, the conclusion is true if
$f\in W^{1,p}_{\rm loc}$ for some $p>n-1$. The latter statement is
the V\"ais\"al\"a result, see Lemma 3 in \cite{Va$_2$}. Theorem
\ref{thOS4.0} is also an extension of the well-known
Gegring-Lehto-Menchoff result in the plane to high dimensions, see,
e.g., \cite{GL}, \cite{LV} and \cite{Me}.
\bigskip

Calderon has shown in \cite{Ca} the preciseness of the condition of
(\ref{eqOS3.001}) for differentiability a.e. of continuous mappings
$f$. Theorem \ref{thOS4.0} shows that we may use the weaker
condition (\ref{eqOS3.0a}) to obtain differentiability a.e. of open
mappings $f$.\bigskip

The condition (\ref{eqOS3.0a}) is not only sufficient but also
necessary for open continuous mappings $f$ $W^{1,\varphi}_{\rm loc}$
from ${\Bbb R}^n$ into ${\Bbb R}^n$, $n\geqslant3$, to have a total
differential a.e. Furthemore, if a function
$\varphi:[0,\infty)\to[0,\infty)$ is increasing, convex and such
that
\begin{equation}\label{eqOS3.0abc}\int\limits_{1}^{\infty}\left[\frac{t}{\varphi(t)}\right]^
{\frac{1}{n-2}}dt=\infty\ ,\end{equation} then there is a
homeomorphism $g:{\Bbb R}^n\to{\Bbb R}^n$, $n\geqslant3$, in the
class $W^{1,\varphi}_{\rm loc}$ which has not a total differential
a.e. Indeed, if $f:{\Bbb R}^{n-1}\to{\Bbb R}$ is a function in the
Calderon construction for $k=n-1$ and
$\varphi_*(t)=\varphi(t+k)-\varphi(k)$, then
\begin{equation}\label{eqOS3.0abcd}\int\limits_{1}^{\infty}\left[\frac{t}{\varphi_*(t)}\right]^
{\frac{1}{n-2}}dt=\infty\end{equation} and $g(x,y)=(x,y+f(x))$,
$x\in{\Bbb R}^{n-1}, y\in{\Bbb R}$, is the desired example because
of $|\nabla g|\leqslant k+|\nabla f|$ and monotonicity of the
function $\varphi$. Thus, the condition (\ref{eqOS3.0a}) already
cannot be weakened even for ho\-meo\-mor\-phisms.
\end{rem}

\cc
\section{The Lusin and Sard properties on surfaces}\label{4}

\begin{theo}{}\label{thOS3.0} {\it Let $\Omega$ be an open set in ${\Bbb R}^k$,
$k\geqslant2$, and let $f:\Omega\to{\Bbb R}^{m},$ $m\ge 1,$ be a
continuous mapping in the class $W^{1,\varphi}(\Omega)$ with an
increasing $\varphi:[0,\infty)\to[0,\infty)$ such that
$\varphi(0)=0$ and
\begin{equation}\label{eqOS3.00}A:=\int\limits_{1}^{\infty}\left[\frac{t}{\varphi(t)}\right]^
{\frac{1}{k-1}}dt<\infty\,.\end{equation} Then
\begin{equation}\label{eqOS3.01}H^{k}(f(E))\leqslant\gamma_{k,m}\,A_*^{k-1}
\int\limits_{E}\varphi_*\left(|\nabla f|\right)\,dm(x)\end{equation}
for every measurable set $E\subset \Omega$ and
$\gamma_{k,m}=(m\alpha_k)^k$ where $\alpha_k$ is a constant from
(\ref{eqOS2.6}) depending only on $k$,
$A_*=A+1/[\varphi(1)]^{1/(k-1)}$, $\varphi_*(0)=0$,
$\varphi_*(t)\equiv\varphi(1)$ for $t\in(0,1)$ and
$\varphi_*(t)=\varphi(t)$ for $t\geqslant1$.} \end{theo}

\bigskip

Thus, we obtain the following conclusions on the Lusin property of
mappings in the Orlicz-Sobolev classes.

\bigskip

\begin{corol}\label{corOS3.00} {\it Under the hypotheses of Theorem \ref{thOS3.0}
the mapping $f$ has the $(N)$-property of Lusin (furthermore, $f$ is
absolutely continu\-ous) with respect to the $k$-dimensional
Hausdorff measure.} \end{corol}

\bigskip

\begin{rem}\label{remOS4.3A}
Note that $H^k({\Bbb R}^m)=0$ for $m<k$ and hence (\ref{eqOS3.01})
is trivial in this case without the condition (\ref{eqOS3.00}).
However, the examples in Section 13 show that the condition
(\ref{eqOS3.00}) is not only sufficient but also necessary for the
(N)--property if $m\ge k$, see Lemma \ref{lemOS13.1} and Remark
\ref{rem13}.
\end{rem}
\bigskip

We obtain also the following consequence of Theorem \ref{thOS3.0} of
the Sard type for mappings in the Orlicz-Sobolev classes, see in
addition Theorem VII.3 in \cite{HW}.
\bigskip

\begin{corol}\label{corOS3.01} {\it Under the hypotheses of Theorem
\ref{thOS3.0}, we have that $H^{k}(f(E))=0$ whenever $|\nabla f|=0$
on a measurable set $E\subset \Omega$ and hence ${\rm dim}_H\
f(E)\leqslant k$} and also ${\rm dim}\ f(E)\leqslant k-1$.
\end{corol}

\bigskip

First such a statement  was established by Sard in \cite{Sar} for
the set of {\bf critical points} of $f$ where $J_f(x)=0$ and then
similar problems studied by many authors for {\bf critical points of
ranks} $r$ where ${\rm rank}\ f'(x)\leqslant r$ and, in particular,
for {\bf supercritical points} where the Jacobian matrix $f'(x)$ is
null at all, see, e.g., \cite{Bat}, \cite{Cr}, \cite{CT$_1$},
\cite{CT$_2$}, \cite{Du}, \cite{Gr}, \cite{Haj}, \cite{Kau},
\cite{No}, \cite{QS}, \cite{Sard$_2$}, \cite{Sard$_3$} and
\cite{Wh}. Usually they requested the corresponding conditions of
smoothness for $f$ without which such statements, generally
speaking, are not true.\medskip

In this connection, we would like to stress here that our result on
supercritical points, Corollary \ref{corOS3.01}, holds without any
assumptions on smoothness of $f$. For instance, this result holds
for all continuous mappings $f$ in the class $W^{1,p}_{\rm loc}$
with $p>k$, see a fine survey on Sard's theorems, in particular, for
Sobolev mappings in the paper \cite{BHS}.

\bigskip

The proof of Theorem \ref{thOS3.0} is based on the following lemma.

\bigskip

\begin{lemma}{}\label{lemOS3.1} {\it Let $\Omega$ be a domain in ${\Bbb R}^k$,
$k\geqslant2$, and let $f:\Omega\to{\Bbb R}^{m}$, $m\ge 1$, be a
continuous mapping in the class $W^{1,\varphi}(G)$ with an
increasing $\varphi:[0,\infty)\to[0,\infty)$ such that
$\varphi(0)=0$ and
\begin{equation}\label{eqOS3.00*}
A:=\int\limits_{1}^{\infty}\left[\frac{t}{\varphi(t)}\right]^
{\frac{1}{k-1}}dt<\infty.\end{equation} Then
\begin{equation}\label{eqOS3.1}{\rm
diam}\,f(C)\leqslant m\alpha_{k}\,A_{*}^{\frac{k-1}{k}}\left[\
\int\limits_{C}\varphi_{*}\left(|\nabla
f|\right)\,dm(x)\right]^{\frac{1}{k}}\end{equation} for every cube
$C\subset \Omega$ whose adges are oriented along coordinate axes
where $\alpha_{k}$ is a constant from (\ref{eqOS2.6}) depending only
on $k$ and $A_*$ and $\varphi_*$ have been defined in Theorem
\ref{thOS3.0}.}
\end{lemma}

\medskip

{\it Proof of Lemma \ref{lemOS3.1}.} Let us prove (\ref{eqOS3.1}) by
induction in $m=1,2,\ldots$. Indeed, (\ref{eqOS3.1}) holds by
Proposition \ref{prOS2.4} and Remark \ref{remOS2.1} for $m=1$ and
$\alpha_k$ from (\ref{eqOS2.6}). Let us assume that (\ref{eqOS3.1})
is valid for some $m=l$ and prove it for $m=l+1$. Consider an
arbitrary vector $\vec{V}=(v_1,v_2,\ldots,v_l,v_{l+1})$ in ${\Bbb
R}^{l+1}$ and the vectors $\vec{V}_1=(v_1,v_2,\ldots,v_l,0)$ and
$\vec{V}_2=(0,\ldots,0,v_{k+1})$. Then
$|\vec{V}|=|\vec{V}_1+\vec{V}_2|\leqslant|\vec{V}_1|+|\vec{V}_2|$.
Thus, denoting by ${\rm Pr}_1\,\vec{V}=\vec{V}_1$ and ${\rm
Pr}_2\,\vec{V}=\vec{V}_2$ the projections of vectors from ${\Bbb
R}^{l+1}$ onto the coordinate hyperplane $y_{l+1}=0$ and on the
$(l+1)$th axis in ${\Bbb R}^{l+1}$, correspondingly, we obtain that
${\rm diam}\, f(C)\leqslant{\rm diam}\,{\rm Pr}_1f(C)+{\rm
diam}\,{\rm Pr}_2f(C)$ and, applying (\ref{eqOS3.1}) under $m=l$ and
$m=1$, we come by monotonicity of $\varphi$ to the inequality
(\ref{eqOS3.1}) under $m=l+1$. The proof is complete.

\medskip

{\it Proof of Theorem \ref{thOS3.0}.} In view of countable
additivity of integral and measure we may assume with no loss of
generality that $E$ is bounded and $\overline{E}\subset G$, i.e.,
$\overline{E}$ is a compactum in $G$. For each $\varepsilon>0$ there
is an open set $\Omega\subset G$ such that $E\subset\Omega$ and
$|\Omega\setminus E|<\varepsilon$, see, e.g., Theorem III (6.6) in
\cite{Sa}. By the above remark we may assume that
$\overline{\Omega}$ is a compactum and, thus, the mapping $f$ is
uniformly continuous in $\Omega$. Hence $\Omega$ can be covered by a
countable collection of closed oriented cubes $C_i$ whose
interiorities are mutually disjoint and such that ${\rm diam}\,
f(C_i)<\delta$ for any prescribed $\delta>0$ and
$\left|\bigcup\limits_{i=1}^{\infty}\partial C_i\right|=0$.\medskip

Thus, by Lemma \ref{lemOS3.1} we have that
$$H_{\delta}^{k}(f(E))\leqslant H_{\delta}^{k}(f(\Omega))\leqslant
\sum\limits_{i=1}^{\infty}\left[\,{\rm diam}\,f(C_i)\right]^{k}
\leqslant$$
$$\leqslant\gamma_{k,m}A_*^{k-1}\int\limits_{\Omega}\varphi_*\left(|\nabla
f|\right)\,dm(x).$$ Finally, by absolute continuity of the
indefinite integral and arbitrariness of $\varepsilon$ and
$\delta>0$ we obtain (\ref{eqOS3.01}).

\bigskip

Combining Proposition \ref{prOS2.1} and Corollary \ref{corOS3.00} we
obtain the follo\-wing statement.

\bigskip

\begin{propo}\label{prOS3.1} {\it Let $k=2,\ldots,n-1$, $U$ be an open set in
${\Bbb R}^n$, $n\geqslant3$, and let $f:U\to{\Bbb R}^m$, $m\ge 1$,
be a continuous mapping in the class $W^{1,\varphi}_{\rm loc}(U)$
for some increasing function $\varphi:[0,\infty)\to[0,\infty)$,
$\varphi(0)=0$, such that (\ref{eqOS3.00}) holds. Then, for every
$k$-dimensional direction $\Gamma$ for a.e. $k$-dimensional plane
${\mathcal P}\in\Gamma$, the restriction of the function $f$ on the
set ${\mathcal P}\cap U$ has the $(N)$-property (furthermore, it is
locally absolutely continuous) with respect to the $k$-dimensional
Hausdorff measure. Moreover, $H^k(f(E))$ $=0$ whenever $\nabla_k
f=0$ on $E\subset P$ for a.e. $P\in\Gamma$.}
\end{propo}

\bigskip

Here $\nabla_k$ denotes the $k$-dimensional gradient of the
restriction of the mapping $f$ to the $k$-dimensional plane $P$.
However, the most important particular case of Proposition
\ref{prOS3.1} for us is the following statement.

\bigskip

\begin{theo}{}\label{thOS3.1} {\it Let $U$ be an open set in ${\Bbb R}^n$,
$n\geqslant3$, and let $\varphi:[0,\infty)\to[0,\infty)$ be an
increasing function such that $\varphi(0)=0$ and
\begin{equation}\label{eqOS3.3}
\int\limits_{1}^{\infty}\left[\frac{t}{\varphi(t)}\right]^
{\frac{1}{n-2}}dt<\infty.\end{equation} Then each continuous mapping
$f:U\to{\Bbb R}^m$, $m\ge 1$, in the class $W^{1,\varphi}_{\rm loc}$
has the $(N)$-property (furthermore, it is locally absolutely
continuous) with respect to the $(n-1)$-dimensional Hausdorff
measure on a.e. hyperplane $\mathcal{P}$ which is parallel to a
fixed hyperplane ${\mathcal P}_0$. Moreover, $H^{n-1}(f(E))=0$
whenever $|\nabla f|=0$ on $E\subset\mathcal{P}$ for a.e. such
$\mathcal{P}$.} \end{theo}

\bigskip

Note that, if the condition (\ref{eqOS3.3}) holds for an increasing
function $\varphi$, then the function $\varphi_*=\varphi(c\,t)$ for
$c>0$ also satisfies (\ref{eqOS3.3}). Moreover, the Hausdorff
measures are quasi-invariant under quasi-isometries. By the
Lindel\"of property of ${\Bbb R}^n$, $U\setminus\{x_0\}$ can be
covered by a countable collection of open segments of spherical
rings in $U\setminus\{x_0\}$ centered at $x_0$ and each such segment
can be mapped onto a rectangular oriented segment of ${\Bbb R}^n$ by
some quasi-isometry, see, e.g., I.5.XI in \cite{Ku} for the
Lindel\"of theorem. Thus, applying Theorem \ref{thOS3.1} piecewise,
we obtain the following conclusion.

\bigskip

\begin{corol}\label{corOS3.2} {\it Under (\ref{eqOS3.3}) each $f\in W^{1,\varphi}_{\rm loc}$
has the $(N)$-property (furthermore, it is locally absolutely
continuous) on a.e. sphere $S$ centered at a prescribed point
$x_0\in{\Bbb R}^n$. Moreover, $H^{n-1}(f(E))=0$ whenever $|\nabla
f|=0$ on $E\subseteq S$ for a.e. such sphere $S$.}
\end{corol}

\bigskip

\begin{rem}\label{remOS3.1} In particular,
(\ref{eqOS3.3}) holds for the functions $\varphi(t)=t^p$, $p>n-1$,
i.e., the given properties hold for $f\in W^{1,p}_{\rm loc}$,
$p>n-1$.\medskip

Note also that (\ref{eqOS3.3}) does not imply the $(N)$-property of
$f:U\to{\Bbb R}^n$ in $U$ with respect to the Lebesgue measure in
${\Bbb R}^n$. The latter conclusion follow, in particular, from the
Ponomarev examples of homeomorphisms $f\in W^{1,p}_{\rm loc}({\Bbb
R}^n)$ for all $p<n$ without $(N)$-property of Lusin, see
\cite{Po}.\medskip

If $m<n-1$, then $H^{n-1}({\Bbb R}^m)=0$ and the (N)--property on
a.e. hyperplane for the mapping $f$ in Theorem \ref{thOS3.1} is
obvious without the condition (\ref{eqOS3.3}). However, the examples
in the final section show that the condition (\ref{eqOS3.3}) are not
only sufficient but also necessary for the $(N)$-property on a.e.
hyperplane if $m\geqslant n-1$, see Remark \ref{rem13} and Theorem
\ref{thOS13.1}.\medskip

The connection of estimates of the Calderon type (\ref{eqOS2.6})
with the (N)--property and differentiability was first found under
the study of the so--called generalized Lipschitzians in the sense
of Rado, see, e.g., \cite{Ca} and V.3.6 in \cite{RR$^*$}, and also
the recent works \cite{BMT}, \cite{KKM}and \cite{RZZ}.
\end{rem}

\cc
\section{On BMO and FMO functions}\label{5}

The BMO space was introduced by John and Nirenberg in \cite{JN} and
soon became one of the main concepts in harmonic analysis, complex
analysis, partial differential equations and related areas, see,
e.g., \cite{HKM} and \cite{ReRy}. \medskip

Let $D$ be a domain in ${\Bbb R}^n$, $n\geqslant 1$. Recall that a
real valued function $\varphi \in$ L$^1_{\mathrm{loc}}(D)$ is said
to be of {\bf boun\-ded mean oscillation} in $D,$ abbr. $\varphi
\in$ BMO(D) or simply $\varphi \in$ {\bf BMO}, if \begin{equation}
\label{eq1.11} \Vert \varphi \Vert _* = \sup\limits_{B\subset D}\ \
\ \dashint\limits_B \vert \varphi (z) - \varphi _B \vert\ dm(z)
\,\,\, < \,\,\, \infty \end{equation} where the supremum is taken
over all balls $B$ in $D$ and \begin{equation}\label{eq1.12} \varphi
_B = \dashint\limits_B \varphi (z)\ dm(z) = {1\over \vert B \vert}
\int\limits_B \varphi (z)\, dm(z) \end{equation} is the mean value
of the function $\varphi  $ over $B.$ Note that L$^{\infty }$(D)
$\subset$ BMO(D) $\subset$ L$^p_{\mathrm{loc}}$(D) for all $1\le
p<\infty ,$ see, e.g., \cite{ReRy}. \medskip

A function $\varphi$ in BMO is said to have {\bf vanishing mean
oscillation}, abbr. $\varphi\in$ {\bf VMO}, if the supremum in
(\ref{eq1.11}) taken over all balls $B$ in $D$ with $|B| <
\varepsilon$ converges to $0$ as $\varepsilon\to 0.$ VMO has been
introduced by Sarason in \cite{Sar}. There are a  number of papers
devoted to the study of partial differential equations with
coefficients of the class { VMO}, see, e.g., \cite{CFL}, \cite{IS},
\cite{MRV$^*$}, \cite{Pal} and \cite{Ra}.\medskip

Following \cite{IR}, we say that a function ${\varphi}: D\to
\mathbb{R} $ has {\bf finite mean oscillation at a point} $z_0 \in
{D}$, write $\varphi\in\ $FMO$(x_0)$, if
\begin{equation}  \label{FMO_eq2.4}
\overline{\lim\limits_{\varepsilon\to 0}}\ \ \ %
\mathchoice {{\setbox0=\hbox{$\displaystyle{\textstyle -}{\int}$}
\vcenter{\hbox{$\textstyle -$}}\kern-.5\wd0}}
{{\setbox0=\hbox{$\textstyle{\scriptstyle -}{\int}$}
\vcenter{\hbox{$\scriptstyle -$}}\kern-.5\wd0}}
{{\setbox0=\hbox{$\scriptstyle{\scriptscriptstyle -}{\int}$}
\vcenter{\hbox{$\scriptscriptstyle -$}}\kern-.5\wd0}}
{{\setbox0=\hbox{$\scriptscriptstyle{\scriptscriptstyle -}{\int}$}
\vcenter{\hbox{$\scriptscriptstyle -$}}\kern-.5\wd0}} \!\int_{B( z_0
,\varepsilon)} |{\varphi}(z)-\tilde{{\varphi}}_{\varepsilon}(z_0)|\
dm(z)\ <\ \infty
\end{equation}
where
\begin{equation}  \label{FMO_eq2.5}
\tilde{{\varphi}}_{\varepsilon}(z_0)=\mathchoice
{{\setbox0=\hbox{$\displaystyle{\textstyle -}{\int}$}
\vcenter{\hbox{$\textstyle -$}}\kern-.5\wd0}}
{{\setbox0=\hbox{$\textstyle{\scriptstyle -}{\int}$}
\vcenter{\hbox{$\scriptstyle -$}}\kern-.5\wd0}}
{{\setbox0=\hbox{$\scriptstyle{\scriptscriptstyle -}{\int}$}
\vcenter{\hbox{$\scriptscriptstyle -$}}\kern-.5\wd0}}
{{\setbox0=\hbox{$\scriptscriptstyle{\scriptscriptstyle -}{\int}$}
\vcenter{\hbox{$\scriptscriptstyle -$}}\kern-.5\wd0}} \!\int_{B( z_0
,\varepsilon)} {\varphi}(z)\ dm(z)
\end{equation}
is the mean value of the function ${\varphi}(z) $ over the ball $B(
z_0
,\varepsilon).$ The condition (\ref{FMO_eq2.4}) includes the assumption that ${%
\varphi}$ is integrable in some neighborhood of the point $z_0.$ We
also say that a function $\varphi$ is of {\bf finite mean
oscillation in the domain} $D$, write $\varphi\in\ $FMO$(D)$ or
simply $\varphi\in\ ${\bf FMO}, if this property holds at every
point $x_0\in D$.
\bigskip

\begin{propo}
\label{FMO_pr2.1} If for some collection of numbers ${\varphi}%
_{\varepsilon}\in {\mathbb{R}},\ \ \varepsilon \in
(0,\varepsilon_0], $
\begin{equation}  \label{FMO_eq2.7}
\overline{\lim\limits_{\varepsilon\to 0}}\ \ \ \mathchoice
{{\setbox0=\hbox{$\displaystyle{\textstyle -}{\int}$}
\vcenter{\hbox{$\textstyle -$}}\kern-.5\wd0}}
{{\setbox0=\hbox{$\textstyle{\scriptstyle -}{\int}$}
\vcenter{\hbox{$\scriptstyle -$}}\kern-.5\wd0}}
{{\setbox0=\hbox{$\scriptstyle{\scriptscriptstyle -}{\int}$}
\vcenter{\hbox{$\scriptscriptstyle -$}}\kern-.5\wd0}}
{{\setbox0=\hbox{$\scriptscriptstyle{\scriptscriptstyle -}{\int}$}
\vcenter{\hbox{$\scriptscriptstyle -$}}\kern-.5\wd0}} \!\int_{B( z_0
,\varepsilon)} |{\varphi}(z)-{\varphi}_{\varepsilon}|\ dm(z)\ <
\infty\, ,
\end{equation}
then $\varphi $ is of finite mean oscillation at $z_0$.
\end{propo}\bigskip

Indeed, by the triangle inequality
$$
\dashint_{B( x_0
,\varepsilon)}\ |\varphi(x)-\overline{\varphi}_{\varepsilon}|\
dm(x)\ \leq \dashint_{B( x_0 ,\varepsilon)}\
|\varphi(x)-\varphi_{\varepsilon}|\ dm(x)\ +\
|\varphi_{\varepsilon}-\overline\varphi_{\varepsilon}|
$$

$$
\leq\
2\cdot\dashint_{B( x_0 ,\varepsilon)}\
|\varphi(x)-\varphi_{\varepsilon}|\ dm(x)\ .
$$

\bigskip

Choosing in Proposition \ref{FMO_pr2.1} $\varphi_{\varepsilon}
\equiv 0, $ $\varepsilon \in (0, \varepsilon_0]$, we obtain the
following statement.\bigskip

\begin{corol} \label{FMO_cor2.1} If for a point $z_0\in{D} $
\begin{equation}  \label{FMO_eq2.8}
\overline{\lim\limits_{\varepsilon\to 0}}\ \ \mathchoice
{{\setbox0=\hbox{$\displaystyle{\textstyle -}{\int}$}
\vcenter{\hbox{$\textstyle -$}}\kern-.5\wd0}}
{{\setbox0=\hbox{$\textstyle{\scriptstyle -}{\int}$}
\vcenter{\hbox{$\scriptstyle -$}}\kern-.5\wd0}}
{{\setbox0=\hbox{$\scriptstyle{\scriptscriptstyle -}{\int}$}
\vcenter{\hbox{$\scriptscriptstyle -$}}\kern-.5\wd0}}
{{\setbox0=\hbox{$\scriptscriptstyle{\scriptscriptstyle -}{\int}$}
\vcenter{\hbox{$\scriptscriptstyle -$}}\kern-.5\wd0}} \!\int_{B( z_0
,\varepsilon)} |{\varphi}(z)|\ dm(z)\ <\ \infty \ ,
\end{equation}
then ${\varphi} $ has finite mean oscillation at $z_0.$
\end{corol}\bigskip

Recall that a point $z_0 \in D$ is called a {Lebesgue point} of a function ${%
\varphi} : D\to\mathbb{R}$ if ${\varphi}$ is integrable in a
neighborhood of $z_0$ and
\begin{equation}  \label{FMO_eq2.7a}
\lim\limits_{\varepsilon\to 0}\ \ \ \mathchoice
{{\setbox0=\hbox{$\displaystyle{\textstyle -}{\int}$}
\vcenter{\hbox{$\textstyle -$}}\kern-.5\wd0}}
{{\setbox0=\hbox{$\textstyle{\scriptstyle -}{\int}$}
\vcenter{\hbox{$\scriptstyle -$}}\kern-.5\wd0}}
{{\setbox0=\hbox{$\scriptstyle{\scriptscriptstyle -}{\int}$}
\vcenter{\hbox{$\scriptscriptstyle -$}}\kern-.5\wd0}}
{{\setbox0=\hbox{$\scriptscriptstyle{\scriptscriptstyle -}{\int}$}
\vcenter{\hbox{$\scriptscriptstyle -$}}\kern-.5\wd0}} \!\int_{B( z_0
,\varepsilon)} |{\varphi}(z)-{\varphi}(z_0)|\ dm(z)\ =\ 0 \ .
\end{equation}
It is known that almost every point in $D$ is a Lebesgue point for
every function ${\varphi}\in L^1 (D).$  We thus have the following
corollary.\bigskip

\begin{corol} \label{FMO_cor2.7b} Every function ${\varphi}:
D\to\mathbb{R},$ which is
locally integrable, has a finite mean oscillation at almost every point in $%
D $. \end{corol}\bigskip

\begin{rem}\rm \label{FMO_rmk2.13a} {\rm Note that the
function ${\varphi }(z)=\log ( 1 / |z| )$ belongs to BMO in the unit disk ${%
\Delta }$, see, e.g., \textrm{\cite{ReRy}, p. 5}, and hence also to FMO. However, $%
\tilde{{\varphi }}_{{\varepsilon }}(0)\rightarrow \infty $ as ${%
\varepsilon }\rightarrow 0,$ showing that the condition
(\ref{FMO_eq2.8}) is only sufficient but not necessary for a
function ${\varphi }$ to be of finite mean oscillation at $z_{0}.$}
\end{rem}\bigskip

Clearly that BMO $\subset$ FMO. By definition FMO$\ \subset
L^1_{\mathrm{loc}}$ but FMO is not a subset of $L^p_{\mathrm{loc}}$
for any $p>1$ in comparison with BMO$_{\mathrm{loc}}\subset
L^p_{\mathrm{loc}}$ for all $p\in [1,\infty)$. Here
BMO$_{\mathrm{loc}}$ stands for the local version of the class BMO.
So, let us give examples showing that FMO is not
BMO$_{\mathrm{loc}}$. \bigskip

{\bf Example 1.} Set $z_n= 2^{-n}, \ r_n= 2^{-pn^2},\ p >1, \
c_n=2^{2n^2},$ $D_n=\{z\in{\Bbb C} : |z-z_n| < r_n \},$ and
$$
\varphi (z)=\sum\limits_{n=1}^{\infty} \ c_n \chi(D_n).
$$
Evidently by Corollary \ref{FMO_cor2.1} that  $\varphi \in FMO({\Bbb
C}\setminus\{0\}). $\medskip

To prove that $\varphi \in FMO(0),$ fix $N$ such that $(p-1)N>1,$
and set $\varepsilon= \varepsilon_N = z_N+ r_N.$ Then
$\bigcup\limits_{n \geq N} \ D_n \ \subset \ {\Bbb D}(\varepsilon)
:= \{z\in{\Bbb C} : |z|<\varepsilon\} $ and
\begin{eqnarray*}
\int\limits_{\Bbb D(\varepsilon)} \varphi\ &=& \ \sum\limits_{n \geq N} \ \int\limits_{D_n} \varphi\ =\ \pi \sum\limits_{n \geq N} c_n r_n^2 \\
&=& \ \sum\limits_{n \geq N} 2^{2(1-p)n^2}\ <\ \sum\limits_{n \geq N} 2^{2(1-p)n}\ \\
&<& \ C \cdot[2^{(1-p)N}]^2\ <\ 2C \varepsilon^2\ .
\end{eqnarray*}
Hence $\varphi \in FMO(0)$ and, consequently, $\varphi \in FMO({\Bbb
C}).$
\par
On the other hand
$$
\int\limits_{\Bbb D(\varepsilon)} \ \varphi^p \ = \ \pi
\sum\limits_{n>N} \ c_n^p \cdot r_n^2 \ = \ \sum\limits_{n>N} 1\ =\
\infty.
$$
Hence $\varphi \notin L^p(\Bbb D(\varepsilon))$ and therefore
$\varphi \notin BMO_{\mathrm{loc}} $ because
BMO$_{\mathrm{loc}}\subset L^p_{\mathrm{loc}}$ for all $p\in
[1,\infty)$.

\bigskip

{\bf Example 2.} We conclude this section by constructing functions
${\varphi} :\mathbb{C} \to \mathbb{R} $ of the class
$C^\infty(\mathbb{C}\backslash\{0\})$ which belongs to FMO but not
to $L^p_{\mathrm{loc}}$ for any $p>1$ and hence not to
BMO$_{\mathrm{loc}}$.\medskip

In this example, $p\ =1+{\delta }$ with an arbitrarily small
${\delta }
>0.$ Set
\begin{equation}  \label{FMO_eq2.11}
{\varphi}_0 (z)= \left \{%
\begin{array}{ll}
e^{\frac {1}{|z|^2-1}}, & \mathrm{if } \ |z| < 1, \\
0, & \mathrm{if} \ |z|\ge 1.%
\end{array}
\right.
\end{equation}

\noindent Then ${\varphi}_0$ belongs to $C^{\infty}_0$ and hence to
$BMO_{\mathrm{loc}}$. Consider the function
\begin{equation}
{\varphi }(z)=\left\{
\begin{array}{ll}
{\varphi }_{k}(z), & \mathrm{if}\ z\in B_{k}, \\
0, & \mathrm{if}\text{ }z\in \mathbb{C}\setminus \bigcup B_{k}%
\end{array}%
\right.  \label{FMO_eq2.11a}
\end{equation}

\noindent where $B_k=B(z_k,r_k), \,\,\, z_k = 2^{-k},\,\,\,
r_k=2^{-(1+{\delta })k^2}, $ ${\delta } > 0,$ and
\begin{equation}  \label{FMO_eq2.12}
{\varphi} _k (z)=2^{2k^2} {\varphi}_0\left(\frac{z-z_k}{r_k}\right),
\ \ z\in B_k,\ \ k=2, 3,\ldots\ .
\end{equation}
Then ${\varphi}$ is smooth in $\mathbb{C}\setminus\{ 0\}$ and,
thus, belongs to BMO$_{\mathrm{loc}}(\mathbb{C}\setminus \{ 0\}),$ and hence to FMO$(%
\mathbb{C}\setminus \{ 0\}).$\medskip

Now,
\begin{equation}  \label{FMO_eq2.5a}
\int\limits_{B_k} {\varphi} _k(z)\ dm(z)\ =\ 2^{-2{\delta } k^2}\int\limits_{%
\mathbb{C}} {\varphi_0} (z)\ dm(z)\ .
\end{equation}
Setting
\begin{equation}  \label{FMO_eq2.13a}
K=K({\varepsilon})= \biggl[\log _2\frac{1}{{\varepsilon}} \biggr]\le \log _2%
\frac{1}{{\varepsilon}}\ ,
\end{equation}
where $[A]$ denotes the integral part of the number $A,$ we have
\begin{equation}  \label{FMO_eq2.13b}
J=\mathchoice {{\setbox0=\hbox{$\displaystyle{\textstyle -}{\int}$}
\vcenter{\hbox{$\textstyle -$}}\kern-.5\wd0}}
{{\setbox0=\hbox{$\textstyle{\scriptstyle -}{\int}$}
\vcenter{\hbox{$\scriptstyle -$}}\kern-.5\wd0}}
{{\setbox0=\hbox{$\scriptstyle{\scriptscriptstyle -}{\int}$}
\vcenter{\hbox{$\scriptscriptstyle -$}}\kern-.5\wd0}}
{{\setbox0=\hbox{$\scriptscriptstyle{\scriptscriptstyle -}{\int}$}
\vcenter{\hbox{$\scriptscriptstyle -$}}\kern-.5\wd0}} \!\int_{D({\varepsilon}%
)} {\varphi}(z)\ dm(z) \le I\cdot {\sum\limits_{k=K}^{\infty}
2^{-2{\delta } k^2}}/{\pi 2^{-2(K+1)}},
\end{equation}
where $I=\int\limits_{\mathbb{C}}{\varphi} (z)\ dm(z) .$ If $
K{\delta } > 1,$ i.e. $K >1/{\delta } ,$ then
\begin{equation}  \label{FMO_eq2.13c}
\sum\limits_{k=K}^{\infty} 2^{-2{\delta } k^2}\le
\sum\limits_{k=K}^{\infty} 2^{-2k} =
2^{-2K}\sum\limits_{k=0}^{\infty} \left( \frac{1}{4}\right)^{k} =
\frac{4}{3}\cdot 2^{-2K},
\end{equation}
i.e., $J\le 16I/3\pi$. Hence
\begin{equation}  \label{FMO_eq2.13}
\overline{\lim\limits_{{\varepsilon}\to 0}}\ \ \ \mathchoice
{{\setbox0=\hbox{$\displaystyle{\textstyle -}{\int}$}
\vcenter{\hbox{$\textstyle -$}}\kern-.5\wd0}}
{{\setbox0=\hbox{$\textstyle{\scriptstyle -}{\int}$}
\vcenter{\hbox{$\scriptstyle -$}}\kern-.5\wd0}}
{{\setbox0=\hbox{$\scriptstyle{\scriptscriptstyle -}{\int}$}
\vcenter{\hbox{$\scriptscriptstyle -$}}\kern-.5\wd0}}
{{\setbox0=\hbox{$\scriptscriptstyle{\scriptscriptstyle -}{\int}$}
\vcenter{\hbox{$\scriptscriptstyle -$}}\kern-.5\wd0}} \!\int_{B({\varepsilon}%
)} {\varphi}(z)\ dm(z)\ <\ \infty\ .
\end{equation}
Thus, ${\varphi}\in$ FMO by Corollary \ref{FMO_cor2.1}.
\par
On the other hand,
\begin{equation}  \label{FMO_eq2.5b}
\int\limits_{B_k} {\varphi} ^{1+{\delta }} _k(z)\ dm(z) = \int\limits_{%
\mathbb{C}} {\varphi}_0 ^{1+{\delta }} (z)\ dm(z)
\end{equation}
and hence ${\varphi}\notin L^{1+{\delta }}(U)$ for any neighborhood
$U$ of $0.$

\cc
\section{On some integral conditions}\label{6}

For every non-decreasing function $\Phi:[0,\infty ]\to [0,\infty ]
,$ the {\bf inverse function} $\Phi^{-1}:[0,\infty ]\to [0,\infty ]$
can be well defined by setting
\begin{equation}\label{eq5.5CC} \Phi^{-1}(\tau)\ =\
\inf\limits_{\Phi(t)\ge \tau}\ t\ .
\end{equation} As usual, here $\inf$ is equal to $\infty$ if the set of
$t\in[0,\infty ]$ such that $\Phi(t)\ge \tau$ is empty. Note that
the function $\Phi^{-1}$ is non-decreasing, too.\medskip

\begin{rem}\label{rmk3.333} Immediately by the
definition  it is evident  that
\begin{equation}\label{eq5.5CCC} \Phi^{-1}(\Phi(t))\ \le\ t\ \ \ \ \
\ \ \ \forall\ t\in[ 0,\infty ]
\end{equation} with the equality in (\ref{eq5.5CCC}) except
intervals of constancy of the function $\Phi(t)$.
\end{rem}

\bigskip

Since the mapping $t\mapsto t^p$ for every positive $p$ is a
sense--preserving homeomorphism $[0, \infty]$ onto $[0, \infty]$ we
may write Theorem 2.1 from \cite{RSY$_1$} in the following form
which is more convenient for further applications. Here, in
(\ref{eq333Y}) and (\ref{eq333F}), we complete the definition of
integrals by $\infty$ if $\Phi_p(t)=\infty ,$ correspondingly,
$H_p(t)=\infty ,$ for all $t\ge T\in[0,\infty) .$ The integral in
(\ref{eq333F}) is understood as the Lebesgue-Stieltjes integral and
the integrals  in (\ref{eq333Y}) and (\ref{eq333B})--(\ref{eq333A})
as the ordinary Lebesgue in\-te\-grals.

\bigskip

\begin{propo} \label{pr4.1aB}{\it\, Let $\Phi:[0,\infty ]\to [0,\infty ]$ be a
non-decreasing function. Set \begin{equation}\label{eq333E} H_p(t)\
=\ \log \Phi_p(t)\ , \qquad \Phi_p(t)=\Phi\left(t^p\right)\,,\quad
p\in (0,\infty)\,.\end{equation}

Then the equality \begin{equation}\label{eq333Y}
\int\limits_{\delta}^{\infty} H^{\,\prime}_p(t)\ \frac{dt}{t}\ =\
\infty
\end{equation} implies the equality \begin{equation}\label{eq333F}
\int\limits_{\delta}^{\infty} \frac{dH_p(t)}{t}\ =\ \infty
\end{equation} and (\ref{eq333F}) is equivalent to
\begin{equation}\label{eq333B} \int\limits_{\delta}^{\infty}H_p(t)\
\frac{dt}{t^2}\ =\ \infty
\end{equation} for some $\delta>0,$ and (\ref{eq333B}) is equivalent to
every of the equalities: \begin{equation}\label{eq333C}
\int\limits_{0}^{\Delta}H_p\left(\frac{1}{t}\right)\ {dt}\ =\ \infty
\end{equation} for some $\Delta>0,$ \begin{equation}\label{eq333D}
\int\limits_{\delta_*}^{\infty} \frac{d\eta}{H_p^{-1}(\eta)}\ =\
\infty
\end{equation} for some $\delta_*>H(+0),$ \begin{equation}\label{eq333A}
\int\limits_{\delta_*}^{\infty}\ \frac{d\tau}{\tau \Phi_p^{-1}(\tau
)}\ =\ \infty \end{equation} for some $\delta_*>\Phi(+0).$
\medskip

Moreover, (\ref{eq333Y}) is equivalent  to (\ref{eq333F}) and hence
(\ref{eq333Y})--(\ref{eq333A})
 are equivalent each to other  if $\Phi$ is in addition absolutely continuous.
In particular, all the conditions (\ref{eq333Y})--(\ref{eq333A}) are
equivalent if $\Phi$ is convex and non--decreasing.}
\end{propo}

\bigskip

It is easy to see that conditions (\ref{eq333Y})--(\ref{eq333A})
become weaker as $p$ increases, see e.g. (\ref{eq333B}). It is
necessary to give one more explanation. From the right hand sides in
the conditions (\ref{eq333Y})--(\ref{eq333A}) we have in mind
$+\infty$. If $\Phi_p(t)=0$ for $t\in[0,t_*]$, then $H_p(t)=-\infty$
for $t\in[0,t_*]$ and we complete the definition $H_p'(t)=0$ for
$t\in[0,t_*]$. Note, the conditions (\ref{eq333F}) and
(\ref{eq333B}) exclude that $t_*$ belongs to the interval of
integrability because in the contrary case the left hand sides in
(\ref{eq333F}) and (\ref{eq333B}) are either equal to $-\infty$ or
indeterminate. Hence we may assume in (\ref{eq333Y})--(\ref{eq333C})
that $\delta>t_0$, correspondingly, $\Delta<1/t_0$ where $t_0\colon
=\sup\limits_{\Phi_p(t)=0}t$, $t_0=0$ if $\Phi_p(0)>0$.\bigskip

Recall that a function  $\Phi :[0,\infty ]\to [0,\infty ]$ is called
{\bf convex} if
$$
\Phi (\lambda t_1 + (1-\lambda) t_2)\ \le\ \lambda\ \Phi (t_1)\ +\
(1-\lambda)\ \Phi (t_2)$$ for all $t_1$ and $t_2\in[0,\infty ] $ and
$\lambda\in [0,1]$.\medskip

Lemma 3.1 from \cite{RSY$_1$} can be written in the following form.

\bigskip

\begin{lemma} \label{lem5.5C} {\it\, Let $Q:{\Bbb B}^n\to [0,\infty ]$ be a measurable
function and let $\Phi:[0,\infty ]\to (0,\infty ]$ be a
non-decreasing convex function.  Then
\begin{equation}\label{eq3.222} \int\limits_{0}^{1}\
\frac{dr}{rq^{\frac{1}{p}}(r)}\ \ge\ \frac{1}{n}\
\int\limits_{eM}^{\infty}\ \frac{d\tau}{\tau \left[\Phi^{-1}(\tau
)\right]^{\frac{1}{p}}}\qquad\qquad\forall\quad p\in (0, \infty)
\end{equation} where $q(r)$ is the average of the function $Q(x)$
over the sphere $|x|=r$ and $M$ is the average of the function
$\Phi\circ Q$ over the unit ball ${\Bbb B}^n$. }\end{lemma}

\bigskip

\begin{rem}\label{rmk3.333A} Note that (\ref{eq3.222}) is
equivalent  for each $p\in (0, \infty)$ to the inequality
\begin{equation}\label{eq3.1!}
\int\limits_0^1\frac{dr}{rq^{\frac{1}{p}}(r)}\ \ge\
\frac{1}{n}\int\limits_{eM}^{\infty}\frac{d\tau}{\tau\Phi_p^{\,-1}(\tau)}\
,\qquad \Phi_p(t)\ \colon =\ \Phi(t^p)\ .
\end{equation}
\end{rem}

\bigskip

\begin{theo} \label{th5.555}{\it\, Let $Q:{\Bbb B}^n\to [0,\infty ]$ be a measurable
function such that \begin{equation}\label{eq5.555}
\int\limits_{{\Bbb B}^n} \Phi (Q(x))\ dm(x)\  <\
\infty\end{equation} where $\Phi:[0,\infty ]\to [0,\infty ]$ is a
non-decreasing convex function such that
\begin{equation}\label{eq3.333a} \int\limits_{\delta_0}^{\infty}\ \frac{d\tau}{\tau
\left[\Phi^{-1}(\tau )\right]^{\frac{1}{p}}}\ =\ \infty\,,\qquad
p\in (0, \infty)\,,
\end{equation} for some $\delta_0\
> \Phi(0).$ Then \begin{equation}\label{eq3.333A}
\int\limits_{0}^{1}\ \frac{dr}{rq^{\frac{1}{p}}(r)}\ =\ \infty
\end{equation} where $q(r)$ is the average of $\ Q(x)$
over the sphere $|x|=r$.}
\end{theo}

\bigskip

\begin{rem}\label{rmk4.7www} In view of Proposition \ref{pr4.1aB},
if (\ref{eq5.555}) holds, then each of the conditions
(\ref{eq333Y})--(\ref{eq333A}) implies the condition
(\ref{eq3.333A}).
\end{rem}

\cc
\section{Moduli of families of surfaces}\label{7}

Let $\omega$ be an open set in $\overline{{\Bbb R}^k}$,
$k=1,\ldots,n-1$. A (continuous) mapping $S:\omega\to{\Bbb R}^n$
is called a $k$-dimensional surface $S$ in ${\Bbb R}^n$. Sometimes
we call the image $S(\omega)\subseteq{\Bbb R}^n$ the surface $S$,
too. The number of preimages
\begin{equation}\label{eq8.2.3}N(S,y)={\rm card}\,S^{-1}(y)={\rm
card}\,\{x\in\omega:\ S(x)=y\},\ y\in{\Bbb R}^n\end{equation} is
said to be a {\bf multiplicity function} of the surface $S$. In
other words, $N(S,y)$ denotes the multiplicity of covering of the
point $y$ by the surface $S$. It is known that the multiplicity
function is lower semicontinuous, i.e.,
$$N(S,y)\ \geqslant\ \liminf_{m\to\infty}\:N(S,y_m)$$
for every sequence $y_m\in{\Bbb R}^n$, $m=1,2,\ldots\,$, such that
$y_m\to y\in{\Bbb R}^n$ as $m\to\infty$; see, e.g., \cite{RR$^*$},
p. 160. Thus, the function $N(S,y)$ is Borel measurable and hence
measurable with respect to every Hausdorff measure $H^k$; see,
e.g., \cite{Sa}, p. 52.

Recall that a $k$-dimensional Hausdorff area in ${\Bbb R}^n$ (or
simply {\bf area}) associated with a surface $S:\omega\to{\Bbb
R}^n$ is given by \begin{equation} \label{eq8.2.4} {\cal
{A}}_S(B)\ =\ {\cal{A}}^{k}_S(B)\ :=\ \int\limits_B N(S,y)\
dH^{k}y \end{equation} for every Borel set $B\subseteq{\Bbb R}^n$
and, more generally, for an arbitrary set that is measurable with
respect to $H^k$ in ${\Bbb R}^n$, cf. 3.2.1 in \cite{Fe}. The
surface $S$ is called {\bf rectifiable} if ${\cal {A}}_S({\Bbb
R}^n)<\infty$, see 9.2 in \cite{MRSY}.

If $\varrho:{\Bbb R}^n\to[0,\infty]$ is a Borel function, then its
{\bf integral over} $S$ is defined by the equality
\begin{equation}\label{eq8.2.5} \int\limits_S \varrho\ d{\cal {A}}\ :=\
\int\limits_{{\Bbb R}^n}\varrho(y)\:N(S,y)\ dH^ky\,.\end{equation}
Given a family $\Gamma$ of $k$-dimensional surfaces $S$, a Borel
function $\varrho:{\Bbb R}^n\to[0,\infty]$ is called {\bf
admissible} for $\Gamma$, abbr. $\varrho\in\mathrm{adm}\,\Gamma$,
if
\begin{equation}\label{eq8.2.6}\int\limits_S\varrho^k\ d{\cal{A}}\ \geqslant\ 1\end{equation}
for every $S\in\Gamma$. Given $p\in(0,\infty)$, the {\bf ${\bf
p}$-modulus} of $\Gamma$ is the quantity
\begin{equation}\label{eq8.2.7}M_p(\Gamma)\ =\
\inf_{\varrho\in\mathrm{adm}\,\Gamma}\int\limits_{{\Bbb
R}^n}\varrho^p(x)\ dm(x).\end{equation} We also set
\begin{equation}\label{eq8.2.8} M(\Gamma)\ =\ M_n(\Gamma)\end{equation} and
call the quantity $M(\Gamma)$ the {\bf modulus of the family}
$\Gamma$. The modulus is itself an outer measure on the collection
of all families $\Gamma$ of $k$-dimensional surfaces.

We say that $\Gamma_2$ is {\bf minorized} by $\Gamma_1$ and write
$\Gamma_2>\Gamma_1$ if every $S\subset\Gamma_2$ has a subsurface
that belongs to $\Gamma_1$. It is known that
$M_p(\Gamma_1)\geqslant M_p(\Gamma_2)$, see \cite{Fu}, p.~176-178.
We also say that a property $P$ holds for {\bf $p$-a.e.} (almost
every) $k$-dimensional surface $S$ in a family $\Gamma$ if a
subfamily of all surfaces of $\Gamma$, for which $P$ fails, has
the $p$-modulus zero. If $0<q<p$, then $P$ also holds for $q$-a.e.
$S$, see Theorem 3 in \cite{Fu}. In the case $p=n$, we write
simply a.e.

\bigskip

\begin{rem}{}\label{rmk8.2.9} The definition of the modulus immediately
implies that, for every $p\in(0,\infty)$ and $k=1,\ldots,n-1$
\medskip

\noindent (1) $p$-a.e. $k$-dimensional surface in ${\Bbb R}^n$ is
rectifiable,
\medskip

\noindent (2) given a Borel set $B$ in ${\Bbb R}^n$ of (Lebesgue)
measure zero, \begin{equation}\label{eq8.2.10} {\cal A}_S(B)\ =\
0\end{equation} for $p$-a.e. $k$-dimensional surface $S$ in ${\Bbb
R}^n$. \end{rem}

\bigskip

The following lemma was first proved in \cite{KR$_1$}, see also
Lemma 9.1 in \cite{MRSY}.

\bigskip

\begin{lemma}{}\label{lem8.2.11} {\it Let $k=1,\ldots,n-1$, $p\in[k,\infty)$,
and let $C$ be an open cube in ${\Bbb R}^n$, $n\geqslant2$, whose
edges are parallel to coordinate axis. If a property $P$ holds for
$p$-a.e. $k$-dimensional surface $S$ in $C$, then $P$ also holds
for a.e. $k$-dimensional plane in $C$ that is parallel to a
$k$-dimensional coordinate plane $H$.} \end{lemma}

\bigskip

The latter a.e. is related to the Lebesgue measure in the
correspon\-ding $(n-k)$-dimensional coordinate plane $H^{\perp}$
that is perpendicular to $H$.

\bigskip

The following statement, see Theorem 2.11 in \cite{KR$_2$} or
Theorem 9.1 in \cite{MRSY}, is an analogue of the Fubini theorem,
cf., e.g., \cite{Sa}, p. 77. It extends Theorem 33.1 in
\cite{Va$_1$}, cf. also Theorem 3 in \cite{Fu}, Lemma 2.13 in
\cite{MRSY$_4$}, and Lemma 8.1 in \cite{MRSY}.

\bigskip

\begin{theo}{}\label{th8.2.12} {\it Let $k=1,\ldots,n-1$, $p\in[k,\infty)$,
and let $E$ be a subset in an open set $\Omega\subset{\Bbb R}^n$,
$n\geqslant2$. Then $E$ is measurable by Lebesgue in ${\Bbb R}^n$
if and only if $E$ is measurable with respect to area on $p$-a.e.
$k$-dimensional surface $S$ in $\Omega$. Moreover, $|E|=0$ if and
only if \begin{equation}\label{eq8.2.13}{\cal A}_S(E)\ =\ 0
\end{equation} on $p$-a.e. $k$-dimensional surface $S$ in
$\Omega$.} \end{theo}

\bigskip

\begin{rem}{}\label{rmk8.2.14} Say by the Lusin theorem, see e.g. Section 2.3.5
in \cite{Fe}, for every measurable function $\varrho:{\Bbb
R}^n\to[0,\infty]$, there is a Borel function $\varrho^*:{\Bbb
R}^n\to[0,\infty]$ such that $\varrho^*=\varrho$ a.e. in ${\Bbb
R}^n$. Thus, by Theorem \ref{th8.2.12}, $\varrho$ is measurable on
$p$-a.e. $k$-dimensional surface $S$ in ${\Bbb R}^n$ for every
$p\in(0,\infty)$ and $k=1,\ldots,n-1$. \end{rem}

\bigskip

We say that a Lebesgue measurable function $\varrho:{\Bbb
R}^n\to[0,\infty]$ is {\bf $p$-extensively admissible} for a
family $\Gamma$ of $k$-dimensional surfaces $S$ in ${\Bbb R}^n$,
abbr. $\varrho\in{\mathrm{ext}}_p\,{\mathrm{adm}}\,\Gamma$, if
\begin{equation}\label{eq8.2.15}\int\limits_S\varrho^k\ d{\cal A}\ \geqslant\ 1
\end{equation} for $p$-a.e. $S\in\Gamma$. The {\bf $p$-extensive
modulus} $\overline M_p(\Gamma)$ of $\Gamma$ is the quantity
\begin{equation}\label{eq8.2.16}\overline M_p(\Gamma)\ =\
\inf\int\limits_{{\Bbb R}^n}\varrho^p(x)\ dm(x)\end{equation}
where the infimum is taken over all
$\varrho\in{\mathrm{ext}}_p\,{\mathrm{adm}}\,\Gamma$. In the case
$p=n$, we use the notations $\overline M(\Gamma)$ and $\varrho\in
{\mathrm{ext}}\,{\mathrm{adm}}\,\Gamma$, respectively. For every
$p\in(0,\infty)$, $k=1,\ldots,n-1$, and every family $\Gamma$ of
$k$-dimensional surfaces in ${\Bbb R}^n$,
\begin{equation}\label{eq8.2.17} \overline M_p(\Gamma)\ =\
M_p(\Gamma)\,.\end{equation}

\cc
\section{Lower and ring $Q$-homeomorphisms}\label{8}

The following concept is motivated by Gehring's ring definition of
quasiconformality in \cite{Ge3}.

Given domains $D$ and $D'$ in $\overline{{\Bbb R}^n}={\Bbb
R}^n\cup\{\infty\}$, $n\geqslant2$,
$x_0\in\overline{D}\setminus\{\infty\}$, and a measurable function
$Q:D\to(0,\infty)$, we say that a homeomorphism $f:D\to D'$ is a
{\bf lower $Q$-homeomorphism at the point} $x_0$ if
\begin{equation}\label{eqOS1.10} M(f\Sigma_{\varepsilon})\
\geqslant\ \inf\limits_{\varrho\in {\mathrm {ext\,
adm}}\,\Sigma_{\varepsilon}}\int\limits_{D\cap
R_{\varepsilon}}\frac{\varrho^n(x)}{Q(x)}\ dm(x)\end{equation} for
every ring $$R_{\varepsilon}\ =\ \{x\in{\Bbb
R}^n:\varepsilon<|\,x-x_0|<\varepsilon_0\}\ ,
\qquad\varepsilon\in(0,\varepsilon_0)\,,\
\varepsilon_0\in(0,d_0)\,,$$ where
\begin{equation}\label{eqOS1.11}d_0\ =\ \sup\limits_{x\in
D}\,|x-x_0|\,,\end{equation} and $\Sigma_{\varepsilon}$ denotes
the family of all intersections of the spheres
$$S(x_0,r)\ =\ \{x\in{\Bbb R}^n:|\,x-x_0|=r\}\ ,\qquad r\in(\varepsilon,\varepsilon_0)\,,$$
with $D$. As usual, the notion can be extended to the case
$x_0=\infty\in\overline{D}$ by applying the inversion $T$ with
respect to the unit sphere in $\overline{{\Bbb R}^n}$,
$T(x)=x/|\,x|^2$, $T(\infty)=0$, $T(0)=\infty$. Namely, a
homeomorphism $f:D\to D'$ is a {\bf lower $Q$-homeo\-mor\-phism
at} $\infty\in\overline{D}$ if $F=f\circ T$ is a lower
$Q_*$-homeomorphism with $Q_*=Q\circ T$ at $0$.

We also say that a homeomorphism $f:D\to{\overline{{\Bbb R}^n}}$
is a {\bf lower $Q$-homeomorphism in} $D$ if $f$ is a lower
$Q$-homeomorphism at every point $x_{0}\in\overline{D}$.

Recall the criterion for homeomorphisms in ${\Bbb R}^n$ to be
lower $Q$-homeo\-morphisms, see Theorem 2.1 in \cite{KR$_1$} or
Theorem 9.2 in \cite{MRSY}.

\medskip

\begin{propo}\label{prOS2.2} {\it Let $D$ and $D'$ be domains in $\overline{{\Bbb R}^n}$,
$n\geqslant2$, let $x_0\in\overline{D}\setminus\{\infty\}$, and
let $Q:D\to(0,\infty)$ a measurable function. A homeomorphism
$f:D\to D'$ is a lower $Q$-homeomorphism at $x_0$ if and only if
\begin{equation}\label{eqOS2.1} M(f\Sigma_{\varepsilon})\ \geqslant\
\int\limits_{\varepsilon}^{\varepsilon_0}
\frac{dr}{||\,Q||\,_{n-1}(r)}\quad\forall\
\varepsilon\in(0,\varepsilon_0)\,,\
\varepsilon_0\in(0,d_0)\,,\end{equation} where
\begin{equation}\label{eqOS2.3} ||Q||\,_{n-1}(r)=\left(\int\limits_{D(x_0,r)}Q^{n-1}(x)\ d{\cal
A}\right)^\frac{1}{n-1}\end{equation} is the $L_{n-1}$-norm of $Q$
over $D(x_0,r)=\{x\in D:|x-x_0|=r\}=D\cap S(x_0,r)$.} \end{propo}

\medskip

Note that the infimum of expression from the right-hand side in
(\ref{eqOS1.10}) is attained for the function
$$\varrho_0(x)=\frac{Q(x)}{||Q||\,_{n-1}(|x|)}\,.$$

\medskip

Now, given a domain $D$ and two sets $E$ and $F$ in
$\overline{{\Bbb R}^n}$, $n\geqslant2$, $\Delta(E,F,D)$ denotes
the family of all paths $\gamma:[a,b]\to\overline{{\Bbb R}^n}$
that join $E$ and $F$ in $D$, i.e., $\gamma(a)\in E$,
$\gamma(b)\in F$, and $\gamma(t)\in D$ for $a<t<b$. Set
\begin{eqnarray} A(r_{1},r_{2},x_{0})&=&\{x\in{\Bbb R}^n:r_{1}<|x-x_{0}|<r_{2}\},
\label{eqOS1.6} \\
S(x_{0},r_{i})&=&\{x\in{\Bbb R}^n:|x-x_{0}|=r_{i}\},\quad i=1,2.
\label{eqOS1.7} \end{eqnarray}

Given domains $D$ in ${\Bbb R}^n$ and $D'$ in $\overline{{\Bbb
R}^n}$, $n\geqslant2$, and a measurable function
$Q:D\to[0,\infty]$, they say that a homeomorphism $f:D\to D'$ is a
{\bf ring $Q$-homeomorphism at a point} $x_{0}\in D$ if
\begin{equation}\label{eqOS1.8} M(\Delta(fS_{1},fS_{2},fD))\ \leqslant
\int\limits_{A}Q(x)\cdot\eta ^{n}(|x-x_{0}|)\ dm(x)\end{equation}
for every ring $A=A(r_{1},r_{2},x_{0})$,
$0<r_{1}<r_{2}<d_{0}={\mathrm{dist}}(x_{0},\partial D)$, and for
every measurable function $\eta:(r_{1},r_{2})\to[0,\infty]$ such
that \begin{equation}\label{eqOS1.9}
\int\limits_{r_{1}}^{r_{2}}\eta(r)\ dr\ =\ 1\,. \end{equation}

The notion was first introduced in the work \cite{RSY$_1$} in the
connection with investigations of the Beltrami equations in the
plane and then it was extended to the space in the work \cite{RS}.

Let us recall the following criterion for ring $Q$-homeomorphisms,
see Theorem 3.15 in \cite{RS} or Theorem 7.2 in \cite{MRSY}.

\medskip

\begin{propo}\label{prOS2.3} {\it Let $D$ be a domain in ${\Bbb R}^n$,
$n\geqslant2$, and $Q:D\to[0,\infty]$ a measurable function. A
homeomorphism $f:D\to{\Bbb R}^n$ is a ring $Q$-homeomorphism at a
point $x_{0}\in D$ if and only if for every
$0<r_{1}<r_{2}<d_{0}={\mathrm{dist}}\,(x_{0},\partial D)$
\begin{equation}\label{eqOS2.4} M(\Delta(fS_{1},fS_{2},fD))\ \leqslant \
\frac{\omega_{n-1}}{I^{n-1}}\end{equation} where $\omega_{n-1}$ is
the area of the unit sphere in ${\Bbb R}^n$, $S_{j}=\{x\in {\Bbb
R}^n:|x-x_{0}|=r_{j}\}$, $j=1,2$, and
$$I\ =\ I(r_{1},r_{2})\ =\ \int\limits_{r_{1}}^{r_{2}}\
\frac{dr}{rq_{x_{0}}^{\frac{1}{n-1}}(r)}$$ and $q_{x_{0}}(r)$ is
the mean value of $Q(x)$ over the sphere $|x-x_{0}|=r$.}
\end{propo}

\medskip

Note that the infimum from the right-hand side in (\ref{eqOS1.8})
holds for the function
$$\eta_{0}(r)=\frac{1}{Irq_{x_{0}}^{\frac{1}{n-1}}(r)}\,.$$

\medskip

\begin{rem}{}\label{rmk2} By the Hesse and Ziemer equalities in \cite{Hes}
and \cite{Zi}, see also the appendixes A3 and A6 in \cite{MRSY},
we have
\begin{equation}\label{eqOS6.1aa} M(\Delta(fS_1,fS_2,fD))\ \leqslant\
\frac{1}{M^{n-1}(f\Sigma)}\end{equation} because
$f\Sigma\subset\Sigma(fS_1,fS_2,fD)$ where $\Sigma$ is a
collection of all spheres centered at $x_0$ between $S_1$ and
$S_2$ and $\Sigma(fS_1,fS_2,fD)$ consists of all
$(n-1)$-dimensional surfaces in $fD$ that separate $fS_1$ and
$fS_2$. \end{rem}

\medskip

Thus, comparing the above criteria for lower and ring
$Q$-ho\-me\-o\-mor\-phisms, we obtain the following conclusion at
inner points.

\medskip

\begin{corol}\label{corOS2.1} {\it Each lower $Q$-homeomorphism $f:D\to D'$ in
${\Bbb R}^n$, $n\geqslant2$, at a point $x_0\in D$ is a ring
$Q^{*}$-homeomorphism with $Q^{*}=Q^{n-1}$ at the point $x_0$.}
\end{corol}

\medskip

\begin{corol}\label{corOS2.2} {\it Each lower $Q$-homeomorphism $f:D\to D'$
in the plane at a point $x_0\in D$ is a ring $Q$-homeomorphism at
the point $x_0$.} \end{corol}

\medskip

It was proved in the work \cite{KPR} that each homeomorphism $f$
of finite distortion in the plane is a lower $Q$-homeomorphism
with $Q(x)=K_{f}(x)$. In the next section we show that the same is
true for a homeomorphism $f$ of finite distortion in ${\Bbb R}^n$,
$n\geqslant3$, if, in addition, $f\in W^{1,\varphi}_{\rm loc}$
where $\varphi$ satisfies the Calderon type condition
(\ref{eqOS3.3}).

\cc
\section{Lower $Q$-homeomorphisms and Orlicz-Sobolev classes}\label{9}

The following statement is key for our further research.

\medskip

\begin{theo}{}\label{thOS4.1} {\it Let $D$ and $D'$ be domains in ${\Bbb R}^n$,
$n\geqslant3$, and let $\varphi:[0,\infty)\to[0,\infty)$ be an
increasing function such that $\varphi(0)=0$ and
\begin{equation}\label{eqOS4.1}
\int\limits_{1}^{\infty}\left[\frac{t}{\varphi(t)}\right]^
{\frac{1}{n-2}}dt<\infty.\end{equation} Then each homeomorphism
$f:D\to D'$ of finite distortion in the class $W^{1,\varphi}_{\rm
loc}$ is a lower $Q$-homeomorphism at every point
$x_0\in\overline{D}$ with $Q(x)=K_f(x)$.} \end{theo}

\bigskip

{\it Proof.} Let $B$ be a (Borel) set of all points $x\in D$ where
$f$ has a total differential $f'(x)$ and $J_f(x)\ne 0$. Then,
applying Kirszbraun's theorem and uniqueness of approximate
differential, see, e.g., 2.10.43 and 3.1.2 in \cite{Fe}, we see that
$B$ is the union of a countable collection of Borel sets $B_l$,
$l=1,2,\ldots\,$, such that $f_l=f|_{B_l}$ are bi-Lipschitz
homeomorphisms, see, e.g., 3.2.2 as well as 3.1.4 and 3.1.8 in
\cite{Fe}. With no loss of generality, we may assume that the $B_l$
are mutually disjoint. Denote also by $B_*$ the rest of all points
$x\in D$ where $f$ has the total differential but with $f'(x)=0$.

By the construction the set $B_0:=D\setminus \left(B\bigcup
B_*\right)$ has Lebesgue measure zero, see Theorem \ref{thOS4.0}.
Hence by Theorem \ref{th8.2.12} ${\cal A}_{S}(B_0)=0$ for a.e.
hypersurface $S$ in ${\Bbb R}^n$ and, in particular, for a.e. sphere
$S_r:=S(x_0,r)$ centered at a prescribed point $x_0\in\overline{D}$.
Thus, by Corollary \ref{corOS3.2} ${\cal A}_{S_r^*}(f(B_0))=0$ as
well as ${\cal A}_{S_r^*}(f(B_*))=0$ for a.e. $S_r$ where
$S_r^*=f(S_r)$.

Let $\Gamma$ be the family of all intersections of the spheres
$S_r$, $r\in(\varepsilon,\varepsilon_0)$,
$\varepsilon_0<d_0=\sup\limits_{x\in D}\,|x-x_0|,$ with the domain
$D.$ Given $\varrho_*\in{\mathrm adm}\,f(\Gamma)$,
$\varrho_*\equiv0$ outside $f(D)$, set $\varrho\equiv 0$ outside
$D$ and on $B_0$
$$\varrho(x)\ \colon=\ \varrho_*(f(x))\Vert f'(x)\Vert \qquad{\rm for}\ x\in D\setminus B_0\,.$$

Arguing piecewise on $B_l$, $l=1,2,\ldots$, we have by 1.7.6 and
3.2.2 in \cite{Fe} that
$$\int\limits_{S_r}\varrho^{n-1}\,d{\cal A}\ \geqslant\
\int\limits_{S_{*}^r}\varrho_{*}^{n-1}\,d{\cal A}\ \geqslant\ 1$$
for a.e. $S_r$ and, thus, $\varrho\in{\mathrm{ext\,adm}}\,\Gamma$.

The change of variables on each $B_l$, $l=1,2,\ldots\,,$ see,
e.g., Theorem 3.2.5 in \cite{Fe}, and countable additivity of
integrals give the estimate
$$\int\limits_{D}\frac{\varrho^n(x)}{K_{f}(x)}\,dm(x)\ \leqslant\
\int\limits_{f(D)}\varrho^n_*(x)\, dm(x)$$ and the proof is
complete.

\bigskip

\begin{corol}{}\label{corOS4.1} {\it Each
homeomorphism $f$ of finite distortion in ${\Bbb R}^n$,
$n\geqslant3$, in the class $W^{1,p}_{\rm loc}$ for $p>n-1$ is a
lower $Q$-homeomorphism at every point $x_0\in\overline{D}$ with
$Q(x)=K_f(x)$.} \end{corol}

\bigskip

\begin{corol}{}\label{corOS4.1*} {\it In particular, each
homeomorphism $f$ of finite distortion in ${\Bbb R}^n$,
$n\geqslant3$, with $K_f\in L^{p}_{\rm loc}$ for $p>n-1$ is a lower
$Q$-homeomorphism at every point $x_0\in\overline{D}$ with
$Q(x)=K_f(x)$.} \end{corol}

\bigskip

\begin{corol}{}\label{corOS4.2} {\it Under the hypotheses of Theorem \ref{thOS4.1},
each homeo\-morphism of finite distortion  $f\in W^{1,\varphi}_{\rm
loc}$, in particular, $f\in W^{1,p}_{\rm loc}$ for $p>n-1$, is a
ring $Q_*$-homeomorphism at every inner point $x_0\in D$ with
$Q_*(x)=\left[K_f(x)\right]^{n-1}$.} \end{corol}

\cc
\section{Equicontinuous and normal families}\label{10}

First of all, recall some general facts on normal families of
mappings in metric spaces. Let $(X,d)$ and $\left(X',d'\right)$ be
metric spaces with distances $d$ and $d'$, respectively. A family
$\frak{F}$ of continuous mappings $f:X\to X'$ is said to be {\bf
normal} if every sequence of mappings $f_m\in\frak{F}$ has a
subsequence $f_{m_k}$ converging uniformly on each compact set
$C\subset X$ to a continuous mapping. Normality is closely related
to the following. A family $\frak{F}$ of mappings $f:X\to X'$  is
said to be {\bf equicontinuous at a point} $x_0\in X$ if for every
$\varepsilon>0$ there is $\delta>0$ such that
$d'(f(x),f(x_0))<\varepsilon$ for all $f\in\frak{F}$ and $x\in X$
with $d(x,x_0)<\delta$. The family $\frak{F}$ is called {\bf
equicontinuous} if $\frak{F}$ is equicontinuous at every point
$x_0\in X$.\bigskip

Given a domain $G$ in ${\Bbb R}^n$, $n\geqslant 2$, and an
increasing function $\varphi:[0,\infty)\to[0,\infty)$  with
$\varphi(0)=0$, $M\in[0,\infty)$ and $x_0\in G$, denote by
$\frak{F}^{\varphi}_M$ the collection of all continuous mappings
$f:G\to{\Bbb R}^m$, $m\geqslant 1$, in the class $W^{1,1}_{\rm loc}$
such that $f(x_0)=0$ and
\begin{equation}\label{eqOS8.111} \int\limits_{G}\varphi\left(|\nabla f|\right)\,dm(x)\leqslant M\ .
\end{equation} By Proposition \ref{prOS2.4} and Remark \ref{remOS2.1} and the
Arzela-Ascoli theorem we obtain the following statement, cf., e.g.,
Theorem 8.1 in \cite{IKO} and Theorem 4.3 in \cite{GMRV}.\bigskip

\begin{corol}\label{corCALDERON} If the function $\varphi$ satisfies
the condition \begin{equation}\label{eq9.111}
\int\limits_{1}^{\infty}\left[\frac{t}{\varphi(t)}\right]^
{\frac{1}{n-1}}dt<\infty\, ,\end{equation} then the class
$\frak{F}^{\varphi}_M$ is equicontinuous and hence normal. If in
addition $\varphi$ is convex, then the class $\frak{F}^{\varphi}_M$
is also closed with respect to the locally uniform convergence.
\end{corol}

\bigskip

Further we give the corresponding theorems for the classes of
homeomorphic mappings under the condition (\ref{eq9.1}) which is
weaker than (\ref{eq9.111}) and without (locally) uniform
constraints of the type (\ref{eqOS8.111}) in these classes.

\bigskip

In what follows, we use in $\overline{{{\Bbb R}}^n}={{\Bbb
R}}^n\bigcup\{\infty\}$ the {\bf spherical (chordal) metric}
$h(x,y)=|\pi(x)-\pi(y)|$ where $\pi$ is the stereographic projection
of $\overline{{{\Bbb R}}^n}$ onto the sphere
$S^n(\frac{1}{2}e_{n+1},\frac{1}{2})$ in ${{\Bbb R}}^{n+1}$:
$$h(x,y)=\frac{|x-y|}{\sqrt{1+{|x|}^2} \sqrt{1+{|y|}^2}}\,,\
x\ne \infty\ne y,\ \ h(x,\infty)=\frac{1}{\sqrt{1+{|x|}^2}}\,.$$
Thus, by definition $h(x,y)\leqslant1$ for all $x$ and
$y\in\overline{{{\Bbb R}}^n}$. The {\bf spherical (chordal)
diameter} of a set $E\subset\overline{{{\Bbb R}}^n}$ is
\begin{equation}\label{eq6.2.16} h(E)=\sup_{x,y\in E}h(x,y)\,.\end{equation}
We use further the following statement of the Arzela-Ascoli type,
see, e.g., Corollary 7.5. in \cite{MRSY}.

\medskip

\begin{propo}\label{prOS9.1} {\it If $(X,d)$ is a separable metric space
and \linebreak $\left(X',d'\right)$ is a compact metric space, then
a family $\frak{F}$ of mappings $f:X\to X'$ is normal if and only if
$\frak{F}$ is equicontinuous.} \end{propo}

\medskip

Combining Theorem \ref{thOS4.1} and Corollaries
\ref{corOS4.1}--\ref{corOS4.2} with the results of the work
\cite{RS}, see also Chapter 7 in \cite{MRSY}, we have the following
statements.

\medskip

\begin{theo}{}\label{th6.4.1} {\it Let $D$ and $D'$ be domains in ${\Bbb
R}^n$, $n\geqslant3$, and let $\varphi:[0,\infty)\to[0,\infty)$ be
an increasing function such that $\varphi(0)=0$ and
\begin{equation}\label{eq9.1} \int\limits_{1}^{\infty}\left[\frac{t}{\varphi(t)}\right]^
{\frac{1}{n-2}}dt<\infty\,.\end{equation} Let $f:D\to D'$ be a
homeomorphism of finite distortion in the Orlicz-Sobolev class
$W^{1,\varphi}_{\rm loc}$ such that $h(\overline{{\Bbb
R}^n}\setminus f(D))\geqslant\Delta>0$. Then, for every $x_0\in D$
and $x\in B(x_0,\varepsilon(x_0))$, $\varepsilon(x_0)<d(x_0)={\rm
dist}(x_0,\partial D)$, \begin{equation}\label{eq6.4.2}
h(f(x),f(x_0))\leqslant\frac{\alpha_n}{\Delta}\,\exp\,
\left\{-\int\limits_{|x-x_0|}^{\varepsilon(x_0)}
\frac{dr}{rk_{x_0}^{\frac{1}{n-1}}(r)}\right\}\end{equation} where
$\alpha_n$ is some constant depending only on $n$ and $k_{x_0}(r)$
is the average of $\ \left[ K_{f}(x)\right]^{n-1}$ over the sphere
$|x-x_0|=r$.}
\end{theo}

\bigskip

\begin{rem}\label{rem6.4.3} The estimate (\ref{eq6.4.2}) can be written in the form
\begin{equation}\label{eq6.4.2a} h(f(x),f(x_0))\leqslant\frac{\alpha_n}{\Delta}\,\exp\,
\left\{-\ \omega^{\frac{1}{n-1}}_{n-1}
\int\limits_{|x-x_0|}^{\varepsilon(x_0)}\frac{dr}{||K_f||_{n-1}(x_0,r)}\right\}\end{equation}
where $||K_f||_{n-1}(x_0,r)$ is the norm of $K_f$ in the space
$L^{n-1}$ over the sphere $|x-x_0|=r$ and $\omega_{n-1}$ is the area
of the unit sphere in ${\Bbb R}^n$. \end{rem}

\medskip

\begin{corol}\label{cor6.4.5a} {\it The estimates (\ref{eq6.4.2})
and (\ref{eq6.4.2a}) hold for homeo\-morphisms $f$ of finite
distortion in the Sobolev classes $W^{1,p}_{\rm loc}$, $p>n-1$. In
particular, these estimates hold for homeomorphisms $f$ of finite
distortion with $K_f\in L^{q}_{\rm loc}$ for $q>n-1$.}
\end{corol}

\medskip

\begin{corol}\label{cor6.4.5} {\it If \begin{equation}\label{eq6.4.6} k_{x_0}(r)\leqslant
\left[\log{\frac{1}{r}}\right]^{n-1}\end{equation} for
$r<\varepsilon(x_0)<\min\left\{1,d(x_0)\right\}$, then
\begin{equation}\label{eq6.4.7} h(f(x),f(x_0))\leqslant\frac{\alpha_n}{\Delta}\,\frac{\log
\frac{1}{\varepsilon(x_0)}}{\log\frac{1}{|x-x_0|}}\end{equation} for
all $x\in B(x_0,\varepsilon(x_0))$.} \end{corol}

\medskip

\begin{corol}\label{cor6.4.8} {\it If \begin{equation}\label{eq6.4.9}
K_f(x)\leqslant\log{\frac{1}{|x-x_0|}},\qquad x\in
B(x_0,\varepsilon(x_0)),\end{equation} then (\ref{eq6.4.7}) holds in
the ball $B(x_0,\varepsilon(x_0))$.} \end{corol}

\medskip

\begin{corol}\label{4.23} {\it Let $n\geqslant3$, $\varphi:[0,\infty)\to[0,\infty)$ be an
increasing function, $\varphi(0)=0$, satisfying (\ref{eq9.1}). Let
$f:{{\Bbb B}}^n\to{\Bbb B}^n$, $f(0)=0$, be a homeomorphism of
finite distortion in the class $W^{1,\varphi}_{\rm loc}$ such that
\begin{equation}\label{eq6.4.24}\int\limits_{\varepsilon<|x|<1}\left[K_f(x)\right]^{n-1}\
\frac{dm(x)}{{|x|}^n}\leqslant c\log{\frac{1}{\varepsilon}}\,,
\qquad \varepsilon\in(0,1).\end{equation} Then
\begin{equation}\label{eq6.4.25}|f(x)|\leqslant\gamma_n\cdot{|x|}^{\beta_n}\end{equation}
where the constants $\gamma_n$ and $\beta_n$ depend only on $n$.}
\end{corol}

\bigskip

\begin{theo}{}\label{theor4} {\it Let $D$ and $D'$ be domains in ${\Bbb R}^n$, $n\geqslant3$, and
let $\varphi:[0,\infty)\to[0,\infty)$ be an increasing function,
$\varphi(0)=0$, such that (\ref{eq9.1}) holds. Suppose $f:D\to D'$
is a homeomorphism of finite distortion in the class
$W^{1,\varphi}_{\rm loc}$ such that $h(\overline{{\Bbb
R}^n}\setminus f(D))\geqslant\Delta>0$ and $K_f(x)\leqslant Q(x)$
where $Q^{n-1}\in{\rm FMO}(x_0)$. Then
\begin{equation}\label{eq54} h(f(x),f(x_0))\leqslant\frac{\alpha_n}{\Delta}
{\left\{{\frac{\log{\frac{1}{\varepsilon_0}}}{\log{\frac{1}{|x-x_0|}}}}\right\}}^{\beta}\quad\forall\
x\in B(x_0,\varepsilon_0)\end{equation} where $\varepsilon_0<{\rm
dist}(x_0,\partial D)$ and $\alpha_n$ depends only on $n$ and
$\beta$ depends on the function $Q$.} \end{theo}

\medskip

\begin{corol}\label{corOS12.*00} {\it In particular, the estimate (\ref{eq54}) holds if
\begin{equation}\label{eqOS12.100} \limsup\limits_{\varepsilon\to0}\dashint_{B(x_0,\varepsilon)}
Q^{n-1}(x)\,dm(x)<\infty\,.\end{equation}} \end{corol}

\medskip

Next, let $D$ be a domain in ${\Bbb R}^n$, $n\geqslant3$, and let
$\varphi:[0,\infty)\to[0,\infty)$ be an increasing function,
$\varphi(0)=0$, $Q:D\to[0,\infty]$ be a measurable function. Let
${\cal O}_{Q,\Delta}^{\varphi}$ be the class of all homeomorphisms
of finite distortion in the Orlicz-Sobolev class $W^{1,\varphi}_{\rm
loc}$ such that $h\left(\overline{{\Bbb R}^n}\setminus f(D)\right)$
$\geqslant\Delta>0$ and $K_f(x)\leqslant Q(x)$ a.e. Moreover, let
${\cal S}^p_{Q,\Delta}$, $p\geqslant1$, denote the classes ${\cal
O}_{Q,\Delta}^{\varphi}$ with $\varphi(t)=t^p$. Finally, let ${\cal
K}^p_{Q,\Delta}$ be the class of all homeomorphisms with finite
distortion such that $K_f\in L^p_{\rm loc}$, $p\geqslant1$,
$K_f(x)\leqslant Q(x)$ a.e. and $h\left(\overline{{\Bbb
R}^n}\setminus f(D)\right)\geqslant\Delta>0$.

\medskip

By Proposition \ref{prOS9.1} the above estimates of distortion now
yield:

\medskip

\begin{theo}{}\label{th6.6.1} {\it Let $\varphi:[0,\infty)\to[0,\infty)$ be an
increasing function such that $\varphi(0)=0$ and (\ref{eq9.1}) hold.
If $Q^{n-1}\in{\rm FMO}$, then ${\cal O}^{\varphi}_{Q,\Delta}$ is a
normal family.} \end{theo}

\medskip

\begin{corol}\label{cor6.6.2} {\it Under (\ref{eq9.1}) the class ${\cal O}^{\varphi}_{Q,\Delta}$
is normal if \begin{equation}\label{eq6.6.3}
\overline{\lim\limits_{\varepsilon\to 0}}\ \
\dashint_{B(x_0,\varepsilon)}Q^{n-1}(x)\,dm(x)<\infty\quad\forall\
x_0\in D.\end{equation}} \end{corol}

\medskip

\begin{corol}\label{cor6.6.4} {\it In particular, the classes ${\cal S}^p_{Q,\Delta}$ and
${\cal K}^p_{Q,\Delta}$ are normal under $p>n-1$ if either
$Q^{n-1}\in {\rm FMO}$ or (\ref{eq6.6.3}) holds.} \end{corol}

\medskip

\begin{theo}{}\label{th6.6.5} {\it Let $\Delta>0$ and $Q:D\to[0,\infty]$
be a measurable function such that \begin{equation}\label{eq6.6.6}
\int\limits_{0}^{\varepsilon(x_0)}\frac{dr}{||Q||_{n-1}(x_0,r)}=\infty\quad\forall\
x_0\in D\end{equation} where $\varepsilon(x_0)<{\rm
dist}(x_0,\partial D)$ and $||Q||_{n-1}(x_0,r)$ denotes the norm of
$Q$ in $L^{n-1}$ over the sphere $|x-x_0|=r$. Then the classes
${\cal O}_{Q,\Delta}^{\varphi}$, ${\cal S}^p_{Q,\Delta}$, ${\cal
K}^p_{Q,\Delta}$ form normal families if $\varphi$ satisfies
(\ref{eq9.1}), correspondingly, $p>n-1$.} \end{theo}

\medskip

\begin{corol}\label{cor6.6.7} {\it The classes ${\cal O}_{Q,\Delta}^{\varphi}$,
${\cal S}^p_{Q,\Delta}$, ${\cal K}^p_{Q,\Delta}$ form normal
fami\-lies if $\varphi$ satisfies (\ref{eq9.1}), correspondingly,
$p>n-1$ and $Q(x)$ has singularities only of the logarithmic type.}
\end{corol}

\medskip

Let $D$ be a fixed domain in the extended space $\overline{{\Bbb
R}^n}={\Bbb R}^n\cup\{\infty\}$, $n\geqslant3$,
$\varphi:[0,\infty)\to[0,\infty)$ be an increasing function,
$\varphi(0)=0$. Given a function $\Phi:[0,\infty]\to[0,\infty]$,
$M>0$, $\Delta>0$, ${\cal O}^{\Phi,\varphi}_{M,\Delta}$ denotes the
collection of all homeomorphisms of finite distortion in the
Orlicz-Sobolev class $W^{1,\varphi}_{\rm loc}$ such that
$h\left(\overline{{\Bbb R}^n}\setminus f(D)\right)\geqslant\Delta>0$
and
\begin{equation}\label{eqOS10.36a} \int\limits_D\Phi\left(K^{n-1}_f(x)\right)
\frac{dm(x)}{\left(1+|x|^2\right)^n}\ \leqslant\ M\,.\end{equation}
Similarly, ${\cal S}^{\Phi,p}_{M,\Delta}$, $p\geqslant1$, denote the
classes ${\cal O}^{\Phi,\varphi}_{M,\Delta}$ with $\varphi(t)=t^p$.
Finally, let ${\cal K}^{\Phi,p}_{M,\Delta}$, $p\geqslant1$, be the
class of all homeomorphisms with finite distortion such that $K_f\in
L^p_{\rm loc}$, $p\geqslant1$, (\ref{eqOS10.36a}) holds for $K_f$
and $h\left(\overline{{\Bbb R}^n}\setminus
f(D)\right)\geqslant\Delta>0$.

\medskip

Combining Theorem \ref{thOS4.1}, Corollaries
\ref{corOS4.1}--\ref{corOS4.2} and also Theorem \ref{th5.555} under
$p=n-1$, we have the following statements, cf. \cite{RS_*}.

\medskip

\begin{theo}{}\label{thOS10.5} {\it Let $\Phi:[0,\infty]\to[0,\infty]$ be a convex
increasing function such that \begin{equation}\label{eqOS10.37a}
\int\limits_{\delta_0}^{\infty}\frac{d\tau}{\tau\left[\Phi^{-1}(\tau)\right]^{\frac{1}{n-1}}}=
\infty\end{equation} for some $\delta_0>\Phi(0)$. Then the classes
${\cal O}^{\Phi,\varphi}_{M,\Delta}$ under (\ref{eq9.1}) and ${\cal
S}^{\Phi,p}_{M,\Delta}$ and ${\cal K}^{\Phi,p}_{M,\Delta}$ under
$p>n-1$ are equicontinuous and, consequently, form normal families
of mappings for every $M\in(0,\infty)$ and $\Delta\in(0,1)$.}
\end{theo}

\medskip

\begin{rem}\label{OS10.5} As it follows from \cite{RS_*}, the
condition (\ref{eqOS10.37a}) is not only sufficient but also
necessary for normality of the given classes. Moreover, by
Proposition \ref{pr4.1aB} we may use instead of (\ref{eqOS10.37a})
each of the equivalent conditions (\ref{eq333Y})--(\ref{eq333A})
under $p=n-1$.\end{rem}

\cc
\section{On domains with regular boundaries}\label{11}

Recall first of all the following topological notion. A domain
$D\subset{\Bbb R}^n$, $n\geqslant2$, is said to be {\bf locally
connected at a point} $x_0\in\partial D$ if, for every
neighborhood $U$ of the point $x_0$, there is a neighborhood
$V\subseteq U$ of $x_0$ such that $V\cap D$ is connected. Note
that every Jordan domain $D$ in ${\Bbb R}^n$ is locally connected
at each point of $\partial D$, see, e.g., \cite{Wi}, p. 66.

\begin{figure}[h]
\centerline{\includegraphics[scale=0.9]{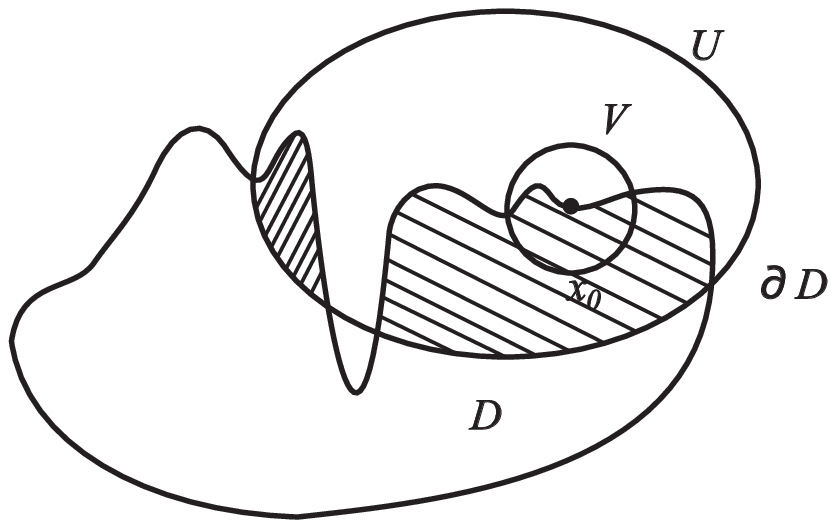}}%
\end{figure}

We say that $\partial D$ is {\bf weakly flat at a point}
$x_0\in\partial D$ if, for every neighborhood $U$ of the point
$x_0$ and every number $P>0$, there is a neighborhood $V\subset U$
of $x_0$ such that
\begin{equation}\label{eq1.5KR}M(\Delta(E,F;D))\geqslant P\end{equation} for
all continua $E$ and $F$ in $D$ intersecting $\partial U$ and
$\partial V$. Here and later on, $\Delta(E,F;D)$ denotes the family
of all paths $\gamma:[a,b]\to\overline{{\Bbb R}^n}$ connecting $E$
and $F$ in $D$, i.e., $\gamma(a)\in E$, $\gamma(b)\in F$ and
$\gamma(t)\in D$ for all $t\in(a,b)$. We say that the boundary
$\partial D$ is {\bf weakly flat} if it is weakly flat at every
point in $\partial D$.

\begin{figure}[h]
\centerline{\includegraphics[scale=0.8]{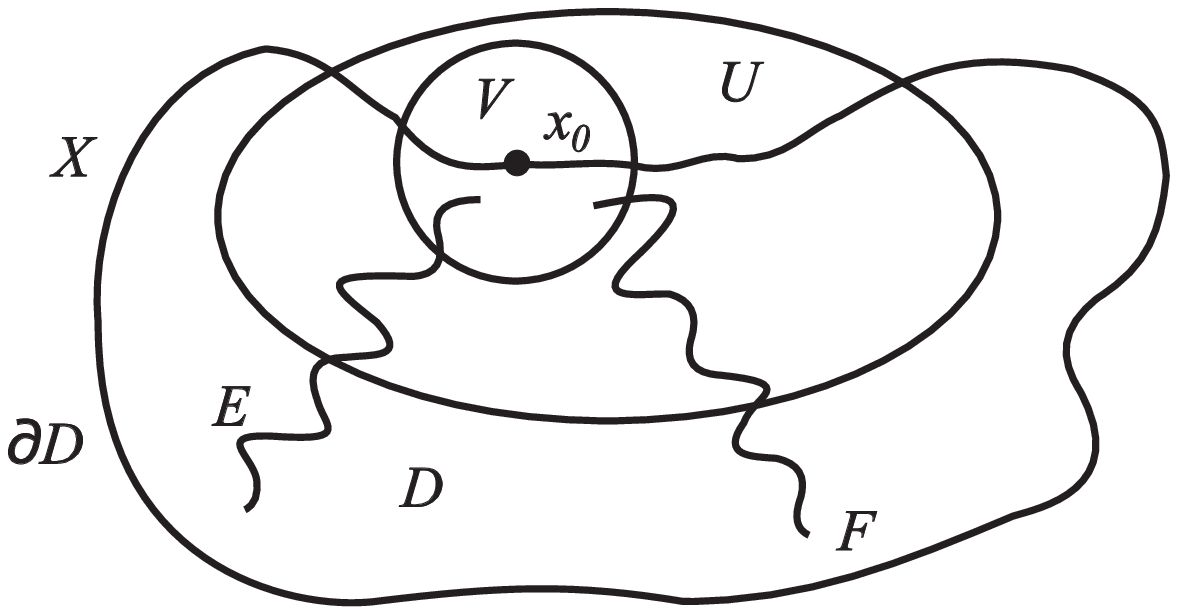}}%
\end{figure}

We also say that a point $x_0\in\partial D$ is {\bf strongly
accessible} if, for every neighborhood $U$ of the point $x_0$,
there exist a compactum $E$ in $D$, a neighborhood $V\subset U$ of
$x_0$ and a number $\delta>0$ such that
\begin{equation}\label{eq1.6KR}M(\Delta(E,F;D))\geqslant\delta\end{equation} for all
continua $F$ in $D$ intersecting $\partial U$ and $\partial V$. We
say that the boundary $\partial D$ is {\bf strongly accessible} if
every point $x_0\in\partial D$ is strongly accessible.

Here, in the definitions of strongly accessible and weakly flat
boun\-daries, one can take as neighborhoods $U$ and $V$ of a point
$x_0$ only balls (closed or open) centered at $x_0$ or only
neighborhoods of $x_0$ in another fundamental system of
neighborhoods of $x_0$. These concepti\-ons can also be extended
in a natural way to the case of $\overline{{\Bbb R}^n}$ and
$x_0=\infty$. Then we must use the corresponding neighborhoods of
$\infty$.

It is easy to see that if a domain $D$ in ${\Bbb R}^n$ is weakly
flat at a point $x_0\in\partial D$, then the point $x_0$ is strongly
accessible from $D$. Moreover, it was proved by us that if a domain
$D$ in ${\Bbb R}^n$ is weakly flat at a point $x_0\in\partial D$,
then $D$ is locally connected at $x_0$, see, e.g., Lemma 5.1 in
\cite{KR$_1$} or Lemma 3.15 in \cite{MRSY}.

The notions of strong accessibility and weak flatness at boundary
points of a domain in ${\Bbb R}^n$ defined in \cite{KR$_0$} are
localizations and generali\-za\-tions of the cor\-res\-pon\-ding
notions introduced in \cite{MRSY$_5$}--\cite{MRSY$_6$}, cf. with
the properties $P_1$ and $P_2$ by V\"ais\"al\"a in \cite{Va$_1$}
and also with the quasicon\-formal accessibility and the
quasiconformal flatness by N\"akki in \cite{Na$_1$}. Many theorems
on a homeomorphic extension to the boundary of quasiconformal
mappings and their generalizations are valid under the condition
of weak flatness of boundaries. The conditi\-on of strong
accessibility plays a similar role for a continuous extensi\-on of
the mappings to the boundary. In particular, recently we have
proved the following significant statements, see either Theorem
10.1 (Lemma 6.1) in \cite{KR$_1$} or Theorem 9.8 (Lemma 9.4) in
\cite{MRSY}.

\bigskip

\begin{propo}\label{prKR2.1} {\it Let $D$ and $D'$ be bounded domains in
${\Bbb R}^n$, $n\geqslant2$, $Q:D\to(0,\infty)$ a measurable
function and $f:D\to D'$ a lower $Q$-homeomorphism on $\partial
D$. Suppose that the domain $D$ is locally connected on $\partial
D$ and that the domain $D'$ has a (strongly accessible) weakly
flat boundary. If
\begin{equation}\label{eqKPR2.11}\int\limits_{0}^{\delta(x_0)} \frac{dr}{||\,Q||_{n-1}(x_0,r)}\
=\ \infty\qquad\forall\ x_0\in\partial D\end{equation} for some
$\delta(x_0)\in(0,d(x_0))$ where $d(x_0)=\sup\limits_{x\in
D}\,|\,x-x_0|$ and
$$||\,Q||_{n-1}(x_0,r)=\left(\int\limits_{D\cap
S(x_0,r)}Q^{n-1}(x)\,d{\cal A}\right)^{\frac{1}{n-1}}\,,$$ then
$f$ has a (continuous) homeomorphic extension $\overline{f}$ to
$\overline{D}$ that maps $\overline{D}$ (into) onto
$\overline{D'}$.}\end{propo}

\medskip

Here as usual $S(x_0,r)$ denotes the sphere $|x-x_0|=r$ and the
closure is understood in the sense of the extended space
$\overline{{\Bbb R}^n}={\Bbb R}^n\cup\{\infty\}$.

\medskip

A domain $D\subset{\Bbb R}^n$ is called a {\bf quasiextremal
distance domain}, abbr. {\bf QED-domain}, see \cite{GM}, if
\begin{equation}\label{e:7.1}M(\Delta(E,F;\overline{{\Bbb R}^n})\leqslant K\cdot
M(\Delta(E,F;D))\end{equation} for some $K\geqslant1$ and all
pairs of nonintersecting continua $E$ and $F$ in $D$.

It is well known, see e.g. Theorem 10.12 in \cite{Va$_1$}, that
\begin{equation}\label{eqKPR2.2}M(\Delta(E,F;{\Bbb R}^n))\geqslant c_n\log{\frac{R}{r}}\end{equation}
for any sets $E$ and $F$ in ${\Bbb R}^n$, $n\geqslant2$,
intersecting all the circles $S(x_0,\rho)$, $\rho\in(r,R)$. Hence
a QED-domain has a weakly flat boundary. One example in
\cite{MRSY}, Section 3.8, shows that the inverse conclusion is not
true even among simply connected plane domains.

A domain $D\subset{\Bbb R}^n$, $n\geqslant2$, is called a {\bf
uniform domain} if each pair of points $x_1$ and $x_2\in D$ can be
joined with a rectifiable curve $\gamma$ in $D$ such that
\begin{equation}\label{e:7.2}s(\gamma)\ \leqslant\ a\cdot|\,x_1-x_2|\end{equation}
and \begin{equation}\label{e:7.3}\min\limits_{i=1,2}\
s(\gamma(x_i,x))\ \leqslant\ b\cdot d(x,\partial D) \end{equation}
for all $x\in\gamma$ where $\gamma(x_i,x)$ is the portion of
$\gamma$ bounded by $x_i$ and $x$, see~\cite{MaSa}. It is known
that every uniform domain is a QED-domain but there exist
QED-domains that are not uniform, see~\cite{GM}. Bounded convex
domains and bounded domains with smooth boundaries are simple
examples of uniform domains and, consequently, QED-doma\-ins as
well as domains with weakly flat boundaries.

A closed set $X\subset{\Bbb R}^n$, $n\geqslant2$, is called a {\bf
null-set for extremal distances}, abbr. {\bf NED-set}, if
\begin{equation}\label{e:8.1} M(\Delta(E,F;{\Bbb R}^n))=M(\Delta(E,F;{\Bbb R}^n\backslash X))\end{equation}
for any two nonintersecting continua $E$ and $F\subset{\Bbb
R}^n\backslash X$.

\medskip

\begin{rem}\label{rmKR2.0} It is known that if $X\subset{\Bbb R}^n$, $n\geqslant2$, is a NED-set, then
\begin{equation}\label{e:8.2}|\,X|=0\end{equation} and $X$ does not locally
disconnect ${\Bbb R}^n$, i.e., see \cite{HW},
\begin{equation}\label{e:8.3}
\dim\,X\leqslant n-2\, ,\end{equation} and, conversely, if a set
$X\subset{\Bbb R}^n$ is closed and
\begin{equation}\label{e:8.4} H^{n-1}(X)=0\,,\end{equation} then $X$
is a NED-set, see \cite{Va$_3$}. Note also that the complement of
a NED-set in ${\Bbb R}^n$ is a very particular case of a
QED-domain. \end{rem}

\bigskip

Further we denote  by $C(X,f)$ the {\bf cluster set} of the mapping
$f:D\to\overline{{\Bbb R}^n}$ for a set $X\subset\overline D$,
\begin{equation}\label{e:8.5}
C(X,f)\ \colon =\left\{y\in\overline{{\Bbb R}^n}:\
y=\lim\limits_{k\to\infty}f(x_k),\ x_k\to x_0\in X,\ x_k\in
D\right\}\ .\end{equation} Note that the inclusion $C(\partial
D,f)\subseteq\partial D'$ holds for every homeomor\-phism $f:D\to
D'$, see, e.g., Proposition 13.5 in \cite{MRSY}.

\cc
\section{The boundary behavior}\label{12}

In this section  we assume that $\varphi:[0,\infty)\to[0,\infty)$ is
an increasing function with $\varphi(0)=0$ such that
\begin{equation}\label{eqOSKRSS}
\int\limits_{1}^{\infty}\left[\frac{t}{\varphi(t)}\right]^
{\frac{1}{n-2}}dt<\infty\ .\end{equation}

In view of Theorem \ref{thOS4.1}, we have by Proposition
\ref{prKR2.1} the following statement.

\bigskip

\begin{theo}\label{thKR9.111} {\it Let $D$ and $D'$ be bounded domains in
${\Bbb R}^n$, $n\geqslant3$, and let $f:D\to D'$ be a homeomorphism
of finite distortion in $W^{1,\varphi}_{\rm loc}$  with the
condition (\ref{eqOSKRSS}). Suppose that the domain $D$ is locally
connected on $\partial D$ and that the domain $D'$ has a (strongly
accessible) weakly flat boundary. If
\begin{equation}\label{eqKPR2.1}\int\limits_{0}^{\delta(x_0)} \frac{dr}{||\,K_f||_{n-1}(x_0,r)}\
=\ \infty\qquad\forall\ x_0\in\partial D\end{equation} for some
$\delta(x_0)\in(0,d(x_0))$ where $d(x_0)=\sup\limits_{x\in
D}\,|\,x-x_0|$ and
$$||\,K_f||_{n-1}(x_0,r)=\left(\int\limits_{D\cap
S(x_0,r)}K_f^{n-1}(x)\,d{\cal A}\right)^{\frac{1}{n-1}}\,,$$ then
$f$ has a (continuous) homeomorphic extension $\overline{f}$ to
$\overline{D}$ that maps $\overline{D}$ (into) onto
$\overline{D'}$.}\end{theo}

\bigskip

In particular, as a consequence of Theorem \ref{thKR9.111} we obtain
the following generalization of the well-known Gehring-Martio
theorem on a homeomorphic extension to the boundary of
quasiconformal mappings between QED domains, cf. \cite{GM}.

\bigskip

\begin{corol}\label{thKPR9.2} {\it Let $D$ and $D'$ be bounded
domains with weakly flat boundaries in ${\Bbb R}^n$, $n\geqslant3$,
and let $f:D\to D'$ be a ho\-meo\-mor\-phism of finite distortion in
$D$ in the class $W^{1,p}_{\rm loc}$, $p>n-1$, in particular,
$K_f\in L^q_{\rm loc}$, $q>n-1$. If the condition (\ref{eqKPR2.1})
holds at every point $x_0\in\partial D$, then $f$ has a homeomorphic
extension to $\overline{D}$.}
\end{corol}

\bigskip

The continuous extenxion to the boundary of the inverse mappings has
a simpler criterion. Namely, in view of Theorem \ref{thOS4.1}, we
have by Theorem 9.1 in \cite{KR$_1$} or Theorem 9.6 in \cite{MRSY}
the next statement.

\bigskip

\begin{theo}\label{thKPR8.2} {\it Let $D$ and $D'$ be domains in ${\Bbb
R}^n$, $n\geqslant3$, $D$ be locally connected on $\partial D$ and
$\partial D'$ be weakly flat. If $f$ is a ho\-meo\-mor\-phism of
finite distortion of $D$ onto $D'$ in the class $W^{1,\varphi}_{\rm
loc}$ with the condition (\ref{eqOSKRSS}) and $K_f\in L^{n-1}(D)$,
then $f^{-1}$ has an extension to $\overline{D'}$ by continuity in
$\overline{{\Bbb R}^n}$.} \end{theo}

\bigskip

However, as it follows from the example in Pro\-po\-si\-tion 6.3
from \cite{MRSY}, any degree of integrability $K_f\in L^q(D)$, $q\in
[1,\infty)$, cannot guarantee  the extension by continuity to the
boundary of the direct mappings.

\bigskip

Similarly, in view of Theorem \ref{thOS4.1}, we have by Theorem 8.1
in \cite{KR$_1$} or Theorem 9.5 in \cite{MRSY} the next result.

\medskip

\begin{theo}\label{thKPR7.1} {\it Let $D$ be a domain in ${\Bbb
R}^n$, $n\ge 3$, $X\subset D$, and let $f$ be a homeomorphism with
finite distortion of $D\backslash X$ into $\overline{{\Bbb R}^n}$ in
$W^{1,\varphi}_{\rm loc}$ with the condition (\ref{eqOSKRSS}).
Suppose that $X$ and $C(X,f)$ are {\rm NED} sets. If
\begin{equation}\label{eq8.9.8} \int\limits_{0}^{\varepsilon({x_0})}
\frac{dr}{||\,K_f||_{n-1}(x_0,r)}\ =\ \infty\ \ \ \ \ \ \ \ \forall\
x_0\in\partial D\end{equation} where $0 < \varepsilon_0 < d_0 = {\rm
dist}\, (x_0,\partial D)$ and
\begin{equation}\label{eq8.9.10}
||\,K_f||_{n-1}(x_0,r)\
=\left(\int\limits_{|x-x_0|=r}K_f^{n-1}(x)\,d{\cal
A}\right)^{\frac{1}{n-1}},\end{equation} then $f$ is extended by
continuity in $\overline{{\Bbb R}^n}$ to $D$.}
\end{theo}

\bigskip

\begin{rem}\label{rmkKPR9.1} In particular, the conclusion of Theorem \ref{thKPR7.1}
is valid if $X$ is a closed set with \end{rem}
\begin{equation}\label{eq8.11.13} H^{n-1}(X)\ =\ 0\ =\
H^{n-1}(C(X,f)).\end{equation}

Finally, in view of Theorem \ref{thOS4.1}, by Theorem
\ref{thKR9.111} as well as by Theorem \ref{th5.555} under $p=n-1$ we
obtain the following result.

\bigskip

\begin{theo}{}\label{thKR4.1} {\it Let $D$ and $D'$ be bounded domains in
${\Bbb R}^n$, $n\geqslant3$, $D$ be locally connected on $\partial
D$ and $D'$ have (strongly accessible) weakly flat boundary. Suppose
$f:D\to D'$ is a homeomorphism of finite distortion in $D$ in the
class $W^{1,\varphi}_{\rm loc}$ with the condition (\ref{eqOSKRSS})
such that
\begin{equation}\label{eqKR4.1}\int\limits_{D}\Phi(K_f^{n-1}(x))\,dm(x)<\infty\end{equation}
for a convex increasing function $\Phi:[0,\infty]\to[0,\infty]$. If
\begin{equation}\label{eqKR4.2}\int\limits_{\delta_0}^{\infty}\frac{d\tau}{\tau\left[\Phi^{-1}(\tau)\right]^{\frac{1}{n-1}}}=
\infty\end{equation} for some $\delta_0>\Phi(0)$, then $f$ has a
(continuous) homeomorphic ex\-ten\-sion $\overline{f}$ to
$\overline{D}$ that maps $\overline{D}$ (into) onto
$\overline{D'}$.}\end{theo}

\bigskip

\begin{rem}\label{rmKR4.1} Note that by Theorem 5.1 and Remark 5.1 in
\cite{KR$_3$} the conditions (\ref{eqKR4.2}) are not only sufficient
but also necessary for continuous extension  to the boundary of $f$
with the integral constraints (\ref{eqKR4.1}).
\medskip

Recall that by Proposition \ref{pr4.1aB}  the condition
(\ref{eqKR4.2}) is equivalent to each of the conditions
(\ref{eq333Y})--(\ref{eq333A}) under $p=n-1$ and, in particular, to
the following condition
\begin{equation}\label{eqKR4.4}\int\limits_{\delta}^{\infty}\log\,\Phi(t)\,\frac{dt}{t^{n'}}\ =\ +\infty\end{equation}
for some $\delta>0$ where $\frac{1}{n'}+\frac{1}{n}=1$, i.e., $n'=2$
for $n=2$, $n'$ is strictly decreasing in $n$ and $n'=n/(n-1)\to1$
as $n\to\infty$.

\bigskip

Finally note that all the results in this section hold, in
particular, if $f\in W^{1,p}_{\rm loc}$, $p>n-1$ and, in particular,
if  $K_f\in L^{q}_{\rm loc}$, $q>n-1$, and, in particular, if $D$
and $D'$ are either bounded convex domains or bounded domains with
smooth boundaries.
\end{rem}

\cc
\section{Some examples}\label{13}

The following lemma is a base for demonstrating preciseness of the
Calderon type conditions in the above results.

\bigskip

\begin{lemma}{}\label{lemOS13.1} {\it Let $\varphi:[0,\infty)\to[0,\infty)$ be
a convex increasing function such that
\begin{equation}\label{eqOS13.1} \int\limits_{1}^{\infty}\left[\frac{t}{\varphi(t)}\right]^
{\frac{1}{k-1}}dt=\infty\end{equation} for a natural number
$k\geqslant2$. Then there is an embedding $g$ of ${\Bbb R}^k$ into
${\Bbb R}^{k+1}$ of the form $g(x)=(x,f(x))$, such that $g\in
W^{1,\varphi}_{\rm loc}$ but $g$ has not $(N)$-property with
respect to $k$-dimensional Hausdorff measure.} \end{lemma}

\bigskip

\begin{corol}\label{corOS13.1} {\it For every $k\geqslant2$, there is an embedding $g$
of ${\Bbb R}^k$ into ${\Bbb R}^{k+1}$ of the form $g(x)=(x,f(x))$
in the class $W^{1,k}_{\rm loc}$ that has not $(N)$-property with
respect to $k$-dimensional Hausdorff measure.} \end{corol}

\bigskip

\begin{rem}\label{rmkOS13.1a} The corresponding examples of embeddings $g$ in the
class $W^{1,k}_{\rm loc}$ for $k=2$ from ${\Bbb R}^{2}$ into ${\Bbb
R}^{3}$ based on the theory of conformal mappings are known long
ago, see, e.g., \cite{Rei} and \cite{Resh}. However, they have not
the form $g(x)=(x,f(x))$ and cannot be applied for constructing
examples of homeomorphisms in the class $W^{1,2}_{\rm loc}$ from
${\Bbb R}^{3}$ into ${\Bbb R}^{3}$ as in Theorem \ref{thOS13.1} and
Corollary \ref{corOS13.2} further. \end{rem}

\bigskip

{\it Proof of Lemma \ref{lemOS13.1}.} We apply further for our
purposes a little modified construction of Calderon in \cite{Ca}, p.
210-211. Let $F$ and $F_*$ be functions from Proposition
\ref{prOS2.5} corresponding to the function
$\varphi_*(t)=\varphi(t+k)-\varphi(k)$. It is clear that $\varphi_*$
satisfies (\ref{eqOS13.1}), too.\medskip

Let us give 3 decreasing sequences of positive numbers $r_l$,
$\varrho_l$ and $\varrho_l^*$, $l=1,2,\ldots$, by induction in the
following way. Set $r_1$ is equal to the maximal number $r>0$ such
that $r\leqslant2^{-2}$ and
$$\int\limits_{|x|\leqslant r}\varphi_*\left(|\nabla
F_*|\right)dm(x)\leqslant2^{-k}.$$ The numbers $\varrho_1$ and
$\varrho_1^*$ are defined from the equalities
$F(\varrho_1)=F(r_1)+1$ and $F(\varrho_1^*)=F(r_1)+3/4$,
correspondingly. If the numbers $r_1,\ldots,r_{l-1}$,
$\varrho_1,\ldots,\varrho_{l-1}$ and
$\varrho_1^*,\ldots,\varrho_{l-1}^*$ have been given, then we set
$r_l$ is equal to the maximal number $r>0$ such that
\begin{equation}\label{eqOS13.3}\int\limits_{|x|\leqslant
r}\varphi_*\left(|\nabla
F_*|\right)dm(x)\leqslant2^{-lk}\end{equation} and, moreover,
\begin{equation}\label{eqOS13.4}r\leqslant
\min\{\, \varrho_{l-1}, 2^{-2l}(\varrho_{l-1}^*-\varrho_{l-1}),
1/[2^{l+2}(l-1)\, |F'(\varrho_{l-1})|]\,\}\,.\end{equation} Then
we define $\varrho_l$ and $\varrho_l^*$ from the equalities
$F(\varrho_l)=F(r_l)+1$ and $F(\varrho_l^*)=F(r_l)+3/4$,
respectively. Note that by monotonicity of the derivative
$|F'(\varrho_{l-1})|=\max\{ |F'(t)|:t\geqslant\varrho_{l-1}\}$. It
is also clear by the construction that $\varrho_l<\varrho_l^*<r_l$
and that the sequence $\varrho_l^*-\varrho_l<r_l$ is decreasing
because the function $F'(t)$ is non-decreasing.

\medskip

Setting $F_l(r)=\min\left[\,1,F(r)-F(r_l)\,\right]$ for
$r\in[0,r_l]$ and $F_l(r)\equiv0$ for $r>r_l$, we see that $F_l$
satisfies a uniform Lipschitz condition, that $F_l(0)=1$ and by
(\ref{eqOS13.3})
\begin{equation}\label{eqOS13.5}\int\limits_{{\Bbb R}^k}\varphi_*\left(|\nabla
F_l^*|\right)dm(x)\leqslant2^{-lk}\end{equation} where
$F_l^*(x)=F_l(|\,x|)$, $x\in{\Bbb R}^k$, $l=1,2,\ldots$.

\medskip

\medskip

Now, denote by $x^l_{j_1,\ldots,j_k}$, $l=1,2,\ldots$,
$j_1,\ldots,j_k=0,\pm1,\pm2,\ldots$, the points in ${\Bbb R}^k$
whose coordinates are integral multiples of $2^{-l}$ with the
natural order in $j_1, ..., j_k$ along the corresponding coordinate
axes. Let $B^l_{j_1,\ldots,j_k}$ be the closed balls centered at
$x^l_{j_1,\ldots,j_k}$ with the radii $r_l$. Note that by the second
condition in (\ref{eqOS13.4}) $r_l\leqslant 2^{-2l}$ and the given
balls are disjoint each to other for every fixed $l=2,\ldots$. Next,
define
$$f_l(x)\ =\ \sum\limits_{j_1,\ldots,j_k}F_l\left(|x-x^l_{j_1,\ldots,j_k}|\right),$$
$$f_p^*(x)\ =\ \sum\limits_{l=1}^{p}2^{-l}f_l(x)$$ and $$f(x)\ =\ \sum\limits_{l=1}^{\infty}2^{-l}f_l(x)\
=\ \lim\limits_{p\to\infty}f_p^*(x)\,.$$ By the construction, the
nonnegative functions $f_l(x)$, $f_p^*(x)$ and $f(x)\leqslant1$ are
continuous. Moreover, it is easy to estimate their oscillations on
the balls $B^p_{j_1,\ldots,j_k}$. In particular,
\begin{equation}\label{eqOS13.13}\osc\limits_{B^p_{j_1,\ldots,j_k}}f^*_{p-1}\leqslant\frac{1}{4}\cdot
2^{-p}\osc\limits_{B^p_{j_1,\ldots,j_k}}f_{p}=2^{-(p+2)}<2^{-(p-1)}.\end{equation}
Indeed, by the triangle inequality and the monotonicity of $F'$,
$$\osc\limits_{B^p_{j_1,\ldots,j_k}}f^*_{p-1}\leqslant
\sum\limits_{l=1}^{p-1}2^{-l}\osc\limits_{B^p_{j_1,\ldots,j_k}}f_{l}\leqslant
r_p\sum\limits_{l=1}^{p-1}|F'(\varrho_l)|\leqslant
r_p(p-1)|F'(\varrho_{p-1})|$$ and, thus, applying the third
condition in (\ref{eqOS13.4}), we come to the
(\ref{eqOS13.13}).\medskip

Let us show that the mapping $g(x)=(x,f(x))$ belongs to the class
$W^{1,\varphi}_{\rm loc}$. To this end, consider an arbitrary closed
oriented unit cube $C$ in ${\Bbb R}^k$ whose vertices have
irrational coordinates. Note that the cube $C$ contains exactly
$2^{lk}$ points $x^l_{j_1,\ldots,j_k}$. Thus, by periodicity of the
picture and the condition (\ref{eqOS13.5}) we have that
\begin{equation}\label{eqOS13.6}\int\limits_{C}\varphi_*\left(|\nabla
f_l|\right)dm(x)\leqslant1\end{equation} and, applying the
(discrete) Jensen inequality, see, e.g., Theorem 86 in \cite{HLP},
we obtain that $$\int\limits_{C}\varphi_*\left(|\nabla
f|\right)dm(x)\leqslant\int\limits_{C}\varphi_*\left(\frac{\sum\limits_{l=1}^{\infty}2^{-l}|\nabla
f_l|}{\sum\limits_{l=1}^{\infty}2^{-l}}\right)dm(x)\leqslant$$
$$\leqslant\sum\limits_{l=1}^{\infty}2^{-l}\int\limits_{C}\varphi_*\left(|\nabla
f_l|\right)dm(x)\leqslant1\,,$$ Finally, since $|\nabla
g|=\sqrt{k+|\nabla f|^2}\leqslant k+|\nabla f|$, we have that
\begin{equation}\label{eqOS13.7}\int\limits_{C}\varphi\left(|\nabla
g|\right)dm(x)\leqslant1+\varphi(k)\,.\end{equation}

Next, let us fix a closed oriented unit cube $C_0$ in ${\Bbb R}^k$
whose center has irrational coordinates and let $E_l$,
$l=1,2,\ldots$, be the union of all balls $B^l_{j_1,\ldots,j_k}$
centered at points $x^l_{j_1,\ldots,j_k}$ in the cube $C_0$. By the
second condition in (\ref{eqOS13.4}) we have that $|E_l|\leqslant
2^{lk}\cdot \Omega_k\cdot 2^{-2lk}=\Omega_k 2^{-lk}$ where
$\Omega_k$ is the volume of the unit ball in ${\Bbb R}^k$. Setting
${\mathcal E}_m=\bigcup\limits_{l=m}^{\infty}E_l$, $m=1,2\ldots$, we
see that $|{\mathcal E}_m|\leqslant \frac{\Omega_k}{2^{k}-1}\
2^{-k(m-1)}\to0$ as $m\to\infty$, i.e., the set ${\mathcal
E}=\bigcap\limits_{m=1}^{\infty}{\mathcal E}_m$ is of the Lebesgue
measure zero in ${\Bbb R}^k$. Similarly, $\mu_{k-1}\left({\rm
pr}_i\, E_l\right)\leqslant\Omega_{k-1}2^{-l(k-1)}$ and
$\mu_{k-1}\left({\rm pr}_i\,{\mathcal E}_m\right)\leqslant
\frac{\Omega_{k-1}}{2^{k-1}-1}\ 2^{-(k-1)(m-1)}\to0$ as
$m\to\infty$, i.e., $\mu_{k-1}({\rm pr}_i\,{\mathcal E})=0$ where
${\rm pr}_i$ denotes the projection into the coordinate hyperplane
$P_i$ which is perpendicular to the $i$-th coordinate axis,
$i=1,2,\ldots,k$ in ${\Bbb R}^k$ and ${\mu}_{k-1}$ is the
$(k-1)$-dimensional Lebesgue measure on $P_i$.

\medskip

Let us prove that every straight line segment $L$ in the cube $C$
which is parallel to a coordinate axis, say to the axis $Ox_1$, and
does not intersect the set ${\mathcal E}$ has only a finite number
of joint points with the spheres $\partial B^{l}_{j_1,\ldots,j_k}$.
Indeed, assume that such a segment $L$ intersects an infinite number
of the closed balls $B^{l}_{j_1,\ldots,j_k}$. Recall that the cube
$C$ intersects only a finite number of such balls under each fixed
$l=1,2,\ldots$. Hence there exists an infinite sequence of balls
$B_{l_m}$ among $B^{l_m}_{j_1,\ldots,j_k}$ such that $L\cap
B_{l_m}\ne\varnothing$, $m=1,2,\ldots$ and $l_{m_1}\ne l_{m_2}$ for
$m_1\ne m_2$, i.e., $l_m\to\infty$ as $m\to\infty$. Note that the
end points of the segment $L$ can belong only to a finite number of
the balls $B_{l_m}$ because in the contrary case it would be
$L\cap{\cal{E}}\ne\varnothing$. Thus, we may assume that
$\mathrm{length}\, L\cap B_{l_m}>0$ for all $m=1,2,\ldots$. Remark
that the distance between the centers $x^{l_m}_{j_1,\ldots,j_k}$ of
the balls $B^{l_m}_{j_1,\ldots,j_k}$, $j_1=0,\pm1,\pm2,\ldots$, as
well as between their projections on the straight line of $L$, is
equal to $2^{-l_m}\to 0$ as $m\to\infty$. Hence we may assume
without loss of generality that the sequence of the segments $L\cap
B_{l_m}$ is monotone decreasing. However, then by the Cantor
theorem, see e.g. 4.41.I ($2'$) in \cite{Ku}, we obtain
$\bigcap\limits_{m=1}^{\infty}B_{l_m}\cap L\ne\varnothing$ that
contradicts the condition ${\cal{E}}\cap L=\varnothing$.

\medskip

Thus, $g$ is piecewise monotone and smooth and hence it is
absolute\-ly continuous on a.e. segment $L$ in the cube $C$ which
is parallel to a coordinate axis. Consequently, $g$ is ${\rm ACL}$
and by (\ref{eqOS13.7}) $g\in W^{1,\varphi}_{\rm loc}$.\bigskip

Let us show that the set $E=g({\mathcal E})$ in ${\Bbb R}^{k+1}$ is
not of $k$-dimensional Hausdorff measure zero. \bigskip

Given a closed oriented unit cube $C_*$ in ${\Bbb R}^k$ whose center
has irrational coordinates and whose edge length $L=2^{-m}$ for some
positive integer $m$, $l\geqslant m$, we have that
\begin{equation}\label{eqOS13.8}\sum\left[d\left(f(B^l_{j_1,\ldots,j_k})\right)\right]^k
\leqslant 2^{3k}L^k \end{equation} where the sum is taken over all
balls $B^l_{j_1,\ldots,j_k}$ whose centers $x^l_{j_1,\ldots,j_k}$
belong to the cube $C_*$. Indeed, $C_*$ contains exactly
$2^{(l-m)k}$ points $x_{j_1,\ldots, j_k}^l$. Moreover, every set
$f(B^l_{j_1,\ldots,j_k})$ is contained in a cylinder whose base
radius is less or equal to $2^{-2l}$, see (\ref{eqOS13.4}), and
whose height is less or equal to
$2^{-(l-1)}+2^{-l}+\ldots=2^{-(l-2)}$, see (\ref{eqOS13.13}). Hence
$$d\left(f(B^l_{j_1,\ldots,j_k})\right)\leqslant\sqrt{2^{-2(l-2)}+2^{2(-2l+1)}}=$$
$$=2^{-(l-2)}\cdot\sqrt{1+2^{-2(l+1)}}\leqslant2^{-(l-3)}=8\cdot2^{-l}$$
that implies (\ref{eqOS13.8}).\medskip

Now, let us prove the following lower estimate of the diameters of
the images of the balls $B^p_{j_1,\ldots,j_k}$ :
\begin{equation}\label{eqOS13.13.13}
d\left(f(B^p_{j_1,\ldots,j_k})\right)\geqslant2^{-(p+1)}.\end{equation}
It is sufficient for this purpose to show that
$$\osc\limits_{B^p_{j_1,\ldots,j_k}}f\geqslant2^{-(p+1)}$$
and, in turn, it suffices to demonstrate that
$$\osc\limits_{L^p_{j_1,\ldots,j_k}}f\geqslant2^{-(p+1)}$$
where $L^p_{j_1,\ldots,j_k}$ is the intersection of the ball
$B^p_{j_1,\ldots,j_k}$ with the straight line $L$ passing through
its center $x^p_{j_1,\ldots,j_k}$ parallely to one of the
coordinate axes. Indeed, by the condition (\ref{eqOS13.4}), the
length of the intersection of the line $L$ with the set ${\cal
{E}}_{p+1}$ can be easy estimated:
$${\rm length}\,(L\cap{\cal{E}}_{p+1})\leqslant\sum\limits_{l=p+1}^{\infty}2r_l\cdot2^l\leqslant2
\sum\limits_{l=p+1}^{\infty}2^{-2l}(\varrho^*_{l-1}-\varrho_{l-1})\cdot2^l\leqslant$$
$$\leqslant2(\varrho^*_{p}-\varrho_{p})\sum\limits_{l=p+1}^{\infty}2^{-l}\leqslant2^{-(p-1)}(\varrho^*_{p}-\varrho_{p})
\leqslant\varrho^*_{p}-\varrho_{p}\,.$$ Hence by the choice of the
number $\varrho^*_{p}$
$$\osc\limits_{\Delta^p_{j_1,\ldots,j_k}}f_p\ \geqslant\
\frac{3}{4}\ \osc\limits_{B^p_{j_1,\ldots,j_k}}f_p$$ where
$\Delta^p_{j_1,\ldots,j_k}=L^p_{j_1,\ldots,j_k}\setminus{\cal
{E}}_{p+1}$. Thus, by the condition (\ref{eqOS13.13}) and the
triangle inequality $$\osc\limits_{\Delta^p_{j_1,\ldots,j_k}}f\ =\
\osc\limits_{\Delta^p_{j_1,\ldots,j_k}}f_p^*\ \geqslant\
\osc\limits_{\Delta^p_{j_1,\ldots,j_k}}f_p-\osc\limits_{\Delta^p_{j_1,\ldots,j_k}}f_{p-1}^*\
\geqslant$$ $$\geqslant\
\frac{1}{2}\cdot2^{-p}\osc\limits_{B^p_{j_1,\ldots,j_k}}f_p\ =\
2^{-(p+1)}$$ and the lower estimate (\ref{eqOS13.13.13}) follows.

\medskip

Finally, let $\varepsilon>0$ and let $\{A_j\}$ be a cover of $E$
such that $d(A_j)<\varepsilon$, $j=1,2,\ldots$. Note that for each
$A_j$ there is a closed oriented cube $C_j$ such that $A_j\subseteq
C_j$ and whose edge length $L_j$ is less or equal to $d(A_j)$.
However, it is more convenient to use closed oriented cubes with
$L_j=2^{-m_j}$ for some positive integer $m_j$ such that
$L_j\leqslant 2d(A_j)$. Let ${\Bbb N}$ be the collection of all
positive integers. Set for arbitrary $l\in{\Bbb N}$
$$S_l=\left\{(j_1,\ldots,j_k):\
x^l_{j_1,\ldots,j_k}\in C_0\right\}\ ,\ \ \ J_l=\left\{j\in{\Bbb N}:
m_j\leqslant l \right\}\ ,$$ and
$$S_l^*=\left\{(j_1,\ldots,j_k)\in S_l:\
x^l_{j_1,\ldots,j_k}\in \bigcup\limits_{j\in J_l}{\rm pr}\,
C_j\right\}\ .$$ Here ${\rm pr}$ denotes the natural projection from
${\Bbb R}^{k+1}$ into ${\Bbb R}^{k}$. Thus, we have by
(\ref{eqOS13.8}) that for every $l\in{\Bbb N}$
$$2^k\sum\limits_{j=1}^{\infty}\left[d(A_j)\right]^k\ \geqslant\
\sum\limits_{j=1}^{\infty}L_j^k\ \geqslant\
2^{-3k}\sum\limits_{(j_1,\ldots,j_k)\in
S_l^*}\left[d\left(f(B^l_{j_1,\ldots,j_k})\right)\right]^k\ .$$
Denote by $N_l$ and $N_l^*$ the numbers of indexes
$(j_1,\ldots,j_k)$ in $S_l$ and $S_l^*$, correspondingly. Note that
by the construction the ratio $N_l^*/N_l$ is non-decreasing and it
converges to $1$ as $l\to\infty$ because $\{C_j\}$ covers $E$,
consequently, $\{{\rm pr}\, C_j\}$ covers $\cal{E}$ and hence
$\bigcup\limits_{j=1}\limits^{\infty}{\rm pr}\, C_j$ includes all
points $x^l_{j_1,\ldots,j_k}$ with $(j_1,\ldots,j_k)\in S_l$.
However, $S_l$ contains exact\-ly $2^{lk}$ indexes
$(j_1,\ldots,j_k)$ and by (\ref{eqOS13.13.13})
$d\left(f(B^l_{j_1,\ldots,j_k})\right)\geqslant2^{-(l+1)}$.
Consequently,
$$\sum\limits_{j=1}^{\infty}\left[d(A_j)\right]^k\geqslant2^{-5k}$$
and, thus, $H^k(E)\geqslant2^{-5k}>0$ in view of arbitrariness of
$\varepsilon>0$. The proof is complete.

\medskip

\begin{rem}\label{rem13} It is known that each homeomorphism of ${\Bbb R}^k$ onto itself in
the class $W^{1,k}_{\rm loc}$ has the $(N)$-property, see Lemma
III.6.1 in \cite{LV} for $k=2$ and \cite{Re$_2$} for $k>2$. The same
is valid also for open mappings, see \cite{MM}. On the other hand,
there exist examples of homeomorphisms $W^{1,p}_{\rm loc}$ for all
$p<k$ that have not the $(N)$-property, see \cite{Po}. Moreover,
Cezari in \cite{Ce} proved that continuous plane mappings
$f:D\to{\Bbb R}^2$ in the class ${\rm ACL}^p$, $p>2$, has the
$(N)$-property and that there exist examples of such mappings in
${\rm ACL}^2$ that have not the $(N)$-property.\medskip

Applying the oblique projection
$h(x)=(x_1,\ldots,x_{k-1},x_k+f(x)/4)$ of the surface
$g(x)=(x,f(x))$, $x\in{\Bbb R}^k$, onto ${\Bbb R}^k$, $k\geqslant2$,
from Lemma \ref{lemOS13.1}, we obtain the corresponding examples of
the continuous mappings $h_k^{\varphi}$ of ${\Bbb R}^k$ onto itself
in the class $W^{1,\varphi}_{\rm loc}$ that have not the
$(N)$-property for convex increasing functions $\varphi$ satisfying
(\ref{eqOS13.1}). In particular, we obtain in this way the example
of a continuous mapping $h_k:{\Bbb R}^k\to{\Bbb R}^k$ in the class
$W^{1,k}_{\rm loc}$ for each integer $k\geqslant2$ without the
$(N)$-property.\medskip

Setting  $H(x,y)=h_{n-1}^{\varphi}(x)$, $x\in{\Bbb R}^{n-1}$,
$y\in{\Bbb R}$, $n\geqslant3$, we obtain examples of continuous
mappings $H:{\Bbb R}^{n}\to{\Bbb R}^{n-1}$ in the class
$W^{1,\varphi}_{\rm loc}$ without the $(N)$-property with respect to
the $(n-1)$-dimensional Hausdorff measure on a.e. hyperplane for
each convex increasing function $\varphi:[0,\infty)\to[0,\infty)$
satisfying the condition (\ref{eqOS13.2}) further. \end{rem}

\begin{theo}{}\label{thOS13.1} {\it Let $\varphi:[0,\infty)\to[0,\infty)$ be
a convex in\-crea\-sing function such that
\begin{equation}\label{eqOS13.2} \int\limits_{1}^{\infty}\left[\frac{t}{\varphi(t)}\right]^
{\frac{1}{n-2}}dt=\infty\end{equation} for a natural number
$n\geqslant3$. Then there is a homeomorphism $H$ of ${\Bbb R}^n$
onto ${\Bbb R}^n$of the form $H(x,y)=(x,y+f(x))$, $x\in{\Bbb
R}^{n-1}$, $y\in{\Bbb R}$, such that $H\in W^{1,\varphi}_{\rm loc}$
but $H$ has not $(N)$-property with respect to $(n-1)$-dimensional
Hausdorff measure on any hyperplane $y={\rm const}$.} \end{theo}

\bigskip

{\it Proof of Theorem \ref{thOS13.1}.} Indeed, the function
$\varphi_*(t):=\varphi(t+1)$ satisfies (\ref{eqOS13.2}). Set
$H(x,y)=g(x)+(0,\ldots,0,y)=(x,y+f(x))$, $x\in{\Bbb R}^{n-1}$,
$y\in{\Bbb R}$, where $g(x)=(x,f(x))$ is the mapping in Lemma
\ref{lemOS13.1} under $k=n-1$ corresponding to the function
$\varphi_*$. Then $|\nabla H|\leqslant1+|\nabla g|$ and by
monotonicity of $\varphi$ we have that $H\in W^{1,\varphi}_{\rm
loc}$ because $g\in W^{1,\varphi_*}_{\rm loc}$.

\bigskip

\begin{corol}\label{corOS13.2} {\it For every $n\geqslant3$, there is a homeomorphism
of ${\Bbb R}^n$ onto ${\Bbb R}^n$ of the class $W^{1,n-1}_{\rm loc}$
of the form $H(x,y)=(x,y+f(x))$, $x\in{\Bbb R}^{n-1}$, $y\in{\Bbb
R}$, without $(N)$-property with respect to $(n-1)$-dimensional
Hausdorff measure on any hyperplane $y={\rm const}$.} \end{corol}

\bigskip

\begin{rem}\label{rmkOS13.1} Note that ${\Bbb R}^n$ can be in the natural way
embedded  into ${\Bbb R}^m$ for each $m>n$. Thus, by Theorems
\ref{thOS3.1} and \ref{thOS13.1} and Remark \ref{rem13}, the
Calderon type condition (\ref{eqOS3.3}) is not only sufficient but
also necessary for continuous mappings $f:{\Bbb R}^n\to{\Bbb R}^m$,
$n\geqslant3$, $m\geqslant n-1$, in the Orlicz-Sobolev classes
$W^{1,\varphi}_{\mathrm{loc}}$ to have the $(N)$-property with
respect to $(n-1)$-dimensional Hausdorff measure on a.e. hyperplane.
Furthermore, Theorem \ref{thOS13.1} shows that the necessity of the
condition (\ref{eqOS3.3}) is valid for $m=n$ even for
ho\-meo\-mor\-phisms $f$. In this connection note also that
Corollaries \ref{corOS13.2} disproves Theorem 1.3 from the preprint
\cite{CHM}.
\end{rem}

\medskip

\noindent Kovtonyuk D., Ryazanov V., Salimov R. and Sevostyanov E. \\
Institute of Applied Mathematics and Mechanics,\\
National Academy of Sciences of Ukraine, \\
74 Roze Luxemburg str., 83114 Donetsk, UKRAINE \\
Phone: +38 -- (062) -- 3110145, Fax: +38 -- (062) -- 3110285 \\
denis$\underline{\ \ }$\,kovtonyuk@bk.ru, vlryazanov1@rambler.ru,\\
salimov@rambler.ru, e$\underline{\ \ }$\,sevostyanov@rambler.ru


\medskip
\medskip

\begin{thebibliography}{99}

\bibitem{Ah} Ahlfors L.: On quasiconformal mappings. J. Anal. Math. \textbf{3},
1--58 (1953/54).

\bibitem{AC} Alberico A. and Cianchi A.: Differentiability properties of
Orlicz-Sobolev functions. Ark. Mat. \textbf{43}, 1--28 (2005).

\bibitem{Am} Ambrosio L.: Metric space valued functions of bounded
variation. Ann. Scuola. Norm. Sup. Pisa Cl. Sci. (4) \textbf{17},
no. 3, 439--478 (1990).

\bibitem{Aseev}
{Aseev V.V.}: {Moduli of families of locally quasisymmetric
 surfaces}. Sib. Math. J. \textbf{30}, no. 3, 353--358 (1989).

\bibitem{AIM} Astala P., Iwaniec T. and Martin G.: Elliptic differential equations
and quasiconformal mappings in the plane. Princeton Math. Ser., vol.
48. Princeton University Press, Princeton (2009).

\bibitem{AIKM} Astala K., Iwaniec T., Koskela P. and Martin G.: Mappings of
BMO-bounded distortion. Math. Ann. \textbf{317}, 703--726 (2000).

\bibitem{BBZ} Bahtin A.K., Bahtina G.P. and Zelinskii Yu.B.:
Topologic-algebraic structures and geometric methods in complex
analysis. Kiev, Inst. Mat. NAHU (2008).

\bibitem{Ba} Balogh Z.M.: Hausdorff dimension distribution of
quasiconformal mappings on the Heisenberg group. J. d'Anal. Math.
\textbf{83}, 289--312 (2001).

\bibitem{BMT} Balogh Z.M., Monti R. and Tyson J.T.: Frequency of
Sobolev and qusiconformal dimension distortion. Research Report
2010-11, 22.07.2010, 1--36 (2010).

\bibitem{Bat} {Bates S.M.}: {On the image size of singular maps}. Proc. AMS
\textbf{114}, 699-705 (1992).

\bibitem{Bel} Belinskii P.P.: General properties of quasiconformal mappings.
Izdat. "Nauka"\ Sibirsk. Otdel., Novosibirsk (1974) [in Russian].

\bibitem{Bi$_1$} Biluta P.A.: Extremal problems for mappings quasiconformal in the
mean. Sib. Mat. Zh. \textbf{6}, 717--726 (1965) [in Russian].


\bibitem{Bi$_2$} Bishop C.J.: Quasiconformal mappings which increase dimension.
Ann. Acad. Sci. Fenn. \textbf{24}, 397--407 (1999).


\bibitem{BO} Birnbaum Z. and Orlicz W.: \"Uber die Verallgemeinerungen des Begriffes der
zueinauder konjugierten Potenzen. Studia. Math. \textbf{3}, 1--67
(1931).

\bibitem{BHS} Bojarski B., Hajlasz P. and Strzelecki P.: Sard's theorem for mappings in H\"older
and Sobolev spaces. Manuscripta Math. \textbf{118}, 383--397 (2005).

\bibitem{Ca} Calderon A.P.: On the differentiability of absolutely
continuous functions. Riv. Math. Univ. Parma \textbf{2}, 203--213
(1951).

\bibitem{Ce} Cesari L.: Sulle transformazioni continue. Annali di
Mat. Pura ed Appl. \textbf{IV 21}, 157--188 (1942).

\bibitem{Ci} Cianchi A.: A sharp embedding theorem for Orlicz-Sobolev spaces.
Indiana Univ. Math. J. \textbf{45} (1), 39--65 (1996).

\bibitem{CFL}
{Chiarenza F., Frasca M. and Longo P.}: {$W^{2,p}$-solvability of
the Dirichlet problem for nondivergence elliptic equations with VMO
coefficients}. Trans. Amer. Math. Soc. \textbf{336}, no. 2, 841--853
(1993).

\bibitem{Cr1} Cristea M.: Open discret mappings having local ${\rm ACL}^n$ inverses.
Preprint, Inst. Math., Romanian Acad., Bucuresti \textbf{6}, 1--29
(2008).

\bibitem{Cr2} Cristea M.: Local homeomorphisms having local ${\rm ACL}^n$ inverses.
Complex Var. Elliptic Equ. \textbf{53} (1), 77--99 (2008).

\bibitem{Cr} {Cristea M.}: {Some properties of the maps in ${\Bbb R}^n$ with applications to the distortion
and singularities of quasiregular mappings}. Rev. Roumaine Math.
Pures Appl. \textbf{36}, 355--368 (1991).

\bibitem{CHM} Cs\"ornyei M., Hencl S. and Maly J.: Homeomorphisms in the Sobolev space
$W^{1,n-1}$. Preprint MATH-KMA-2007/252, Prague Charles Univ., 1--15
(2007).

\bibitem{CT$_1$} {Church P. T. and Timourian J. G.}: {Differentiable Maps with Small Critical Set or
Critical Set Image}. Indiana Univ. Math. J. \textbf{27}, 953--971
(1978).


\bibitem{CT$_2$} {Church P. T. and Timourian J. G.}: {Maps having
0--dimensional critical set image}. Indiana Univ. Math. J. ,
\textbf{27}, 813--832 (1978).

\bibitem{Do} Donaldson T.: Nonlinear elliptic boundary-value problems in
Orlicz-Sobolev spaces. J. Diff. Eq. \textbf{10}, 507--528 (1971).

\bibitem{Dub} Dubinin V.N.: Capacities of condencers  and symmetrization
in geometric function theory of complex variable. Vladivostok,
Dal'nauka (2009).

\bibitem{Du} Dubovickii A.Yu.: On the structure of level sets of
differentiable mappings of an $n-$dimensional cube into a
$k-$dimensional cube. Izv. Akad. Nauk SSSR. Ser. Mat. \textbf{21},
371--408 (1957).

\bibitem{Fa} Fadell A.G.: A note on a theorem of Gehring and
Lehto. Proc. Amer. Math. Soc. \textbf{49}, 195--198 (1975).

\bibitem{FKZ} Faraco D., Koskela P. and Zhong X.: Mappings of finite distortion:
The degree of regularity. Adv. Math. \textbf{190} (2), 300--318
(2005).

\bibitem{Fe} Federer H.: Geometric Measure Theory. Springer-Verlag, Berlin
(1969).

\bibitem{Fu} Fuglede B.: Extremal length and functional completion. Acta Math.
\textbf{98}, 171--219 (1957).

\bibitem{Ge$_2$} Gehring F.W.: Symmetrization of rings in space. Trans. Amer. Math.
Soc. \textbf{101}, 499--519 (1961).

\bibitem{Ge3} Gehring F.W.: Rings and quasiconformal mappings in space. Trans.
Amer. Math. Soc. \textbf{103}, 353--393 (1962).

\bibitem{GI} Gehring F.W. and Iwaniec T.: The limit of mappings with finite
distortion, Ann. Acad. Sci. Fenn. Math. \textbf{24}, 253--264
(1999).

\bibitem{GL} Gehring F.W. and Lehto O.: On the total
differentiability of functions of a complex variable. Ann. Acad.
Sci. Fenn. Ser. A1. Math. \textbf{272}, 3--8 (1959).

\bibitem{GM} Gehring F.W. and Martio O.: Quasiextremal distance domains and
extension of quasiconformal mappings. J. Anal. Math. \textbf{45},
181--206 (1985).

\bibitem{GV} Gehring F.W. and V\"ais\"al\"a J.: Hausdorff dimension and quasiconformal mappings.
J. London Math. Soc. (2) \textbf{6}, 504--512 (1973).

\bibitem{Go}
{Golberg A.}: {Homeomorphisms with finite mean dilatations}.
Contemporary Math. \textbf{382}, 177--186 (2005).

\bibitem{GK} Golberg A. and Kud'yavin V.S.: Mean coefficients of quasiconformality
of pair of domains. Ukrain. Mat. Zh. \textbf{43} (12), 1709--1712
(1991) [in Russian]; translation in Ukrain. Math. J. \textbf{43},
1594--1597 (1991).

\bibitem{GR} Gol’dshtein V. M. and  Reshetnyak Yu.G.: Introduction to the theory
of functions with distributional derivatives and quasiconformal
mappings: Nauka, Moscow (1983); English transl., Quasiconformal
mappings and Sobolev spaces. Kluwer, Dordrecht (1990).

\bibitem{GM$_*$} Gossez J.-P. and Mustonen V.: Variational inequalities in
Orlicz-Sobolev spaces. Nonlinear Anal. Theory Meth. Appl.
\textbf{11}, 379--392 (1987).

\bibitem{Gr} Grinberg E.L.: On the smoothness hypothesis in Sard's
theorem. Amer. Math. Monthly \textbf{92}, no. 10, 733--734 (1985).



\bibitem{GMSV}
{Gutlyanskii V., Martio O., Sugawa T. and Vuorinen M.}: {On the
degenerate Beltrami equation}. Trans. Amer. Math. Soc. \textbf{357},
875--900 (2005).

\bibitem{GMRV}
{Gutlyanski\u{i} V.Ya., Martio O., Ryazanov V.I. and Vuorinen M.}:
{On convergence  theorems for space quasiregular mappings}. Forum
Math. \textbf{10}, 353--375 (1998).

\bibitem{GRSY}
{Gutlyanskii V.,Ryazanov V., Srebro U. and Yakubov E.}: {On recent
advances in the Beltrami equations}. Ukrainian Math. Bull.
\textbf{7}, no. 4, 875--900 (2005).

\bibitem{Ha} {Hajlasz P.}: {Sobolev spaces on an arbitrary metric space}.
Potential Anal. \textbf{5}, 403--415 (1996).

\bibitem{Haj} {Hajlasz P.}: {Whitney's example by way of Assouad's embedding}.
Proc. Amer. Math. Soc. \textbf{131}, 3463--3467 (2003).

\bibitem{HLP} Hardy G.H., Littlewood J.E. and Polia G:
Inequalities. Cambridge University Press, Cambridge (1934).

\bibitem{He} {Heinonen J.}: {Lectures of Analysis on Metric Spaces}. Springer,
New York etc. (2000).

\bibitem{HKM} {Heinonen J., Kilpelainen T. and Martio O.}: { Nonlinear Potential
Theory of Degenerate Elliptic Equations}. Oxford Mathematical
Monographs. Clarendon Press, Oxford-New York-Tokio (1993).

\bibitem{HK$^*_2$} {Heinonen J. and Koskela P.}: {Quasiconformal maps in metric spaces
with controlled geometry} Acta Math. \textbf{181}, 1--41 (1998).

\bibitem{HK$^*_1$} Heinonen J. and Koskela P.: Sobolev mappings with integrable
dilatations. Arch. Rational Mech. Anal. \textbf{125}, 81--97
(1993).

\bibitem{HKST} Heinonen J., Koskela P., Shanmugalingam N. and Tyson J.T.:
Sobolev spaces of Banach space--valued functions and quasiconformal
mappings. J. Anal. Math. \textbf{85}, 87--139 (2001).

\bibitem{HKO} Hencl S., Koskela P. and Onninen J.: A note on extremal mappings of
finite distortion. Math. Res. Lett. \textbf{12} (2-3), 231--237
(2005).

\bibitem{HM$^*$} Hencl S. and Maly J.: Mappings with finite distortion: Hausdorff
measure of zero sets. Math. Ann. \textbf{324} (3), 451--464
(2002).

\bibitem{HK$_1$} Herron D.A. and Koskela P.: Locally uniform domains and
quasiconformal mappings. Ann. Acad. Sci. Fenn. Ser. A1. Math.
\textbf{20}, 187--206 (1995).

\bibitem{Hes} Hesse J.: A $p$-extremal length and $p$-capacity equality. Ark.
Mat. \textbf{13}, 131--144 (1975).

\bibitem{HP} Holopainen I. and Pankka P.: Mappings of finite distortion: Global
homeomorphism theorem. Ann. Acad. Sci. Fenn. Math. \textbf{29}(1),
59--80 (2004).

\bibitem{Hs} Hsini M.: Existence of solutions to a semilinear elliptic system through
generalized Orlicz-Sobolev spaces. J. Partial Differ. Equ.
\textbf{23} (2), 168--193 (2010).

\bibitem{HW} Hurewicz W. and Wallman H.: Dimension theory. Princeton Univ. Press,
Princeton, NJ (1948).

\bibitem{IR} Ignat'ev A. and Ryazanov V.: Finite mean oscillation in the mapping theory.
Ukrainian Math. Bull. \textbf{2} (3), 403--424 (2005).

\bibitem{IKO$_1$} Iwaniec T., Koskela P. and Onninen J.: Mappings of finite distortion:
Monotonicity and continuity. Invent. Math. \textbf{144} (3),
507--531 (2001).

\bibitem{IKO} Iwaniec T., Koskela P. and Onninen J.: Mappings of
finite distortion: Compactness. Ann. Acad. Sci. Fenn. Math.
\textbf{27} (2), 391--417 (2002).

\bibitem{IM} Iwaniec T. and Martin G.: Geometric function theory and non-linear analysis.
Oxford Math. Monogr., Oxford Univ. Press, Oxford (2001).

\bibitem{IS} Iwaniec T. and Sverak V.: On mappings with integrable
dilatation. Proc. Amer. Math. Soc. \textbf{118}, 181--188 (1993).

\bibitem{JN} John F. and Nirenberg L.: On functions of
bounded mean oscillation. Comm. Pure Appl. Math. \textbf{14},
415--426 (1961).

\bibitem{Ka} Kallunki S.: Mappings of finite distortion: The metric definition.
Dissertation, Univ. Jyvaskyla, Jyvaskyla. Ann. Acad. Sci. Fenn.
Math. Diss. \textbf{131}, pp. 33 (2002).

\bibitem{Kau} Kaufman R.: A singular map of a cube onto a square. J. Diff. Geom.
\textbf{14}, 593--594 (1979).

\bibitem{KKM} Kauhanen J., Koskela P. and Maly J.: On functions with derivatives in a Lorentz space.
Manuscripta math. \textbf{10}, 87--101 (1999).

\bibitem{KKM$_1$} Kauhanen J., Koskela P. and Maly J.: Mappings of finite distortion:
Discreteness and openness. Arch. Rat. Mech. Anal. \textbf{160},
135--151 (2001).

\bibitem{KKM$_2$} Kauhanen J., Koskela P. and Maly J.: Mappings of finite distortion:
Condition $N$. Michigan Math. J. \textbf{49}, 169--181 (2001).

\bibitem{KKMOZ} Kauhanen J., Koskela P., Maly J., Onninen J. and Zhong X.: Mappings
of finite distortion: Sharp Orlicz-conditions. Rev. Mat.
Iberoamericana \textbf{19} (3), 857--87 (2003).

\bibitem{KP} Khruslov E.Ya. and Pankratov L.S.: Homogenization of the Dirichlet variational
problems in Sobolev-Orlicz spaces. Operator theory and its
applications (Winuipeg, MB, 1998), 345-366, Fields Inst. Commun.,
25, Amer. Math. Soc,. Providence, RI, 2000.

\bibitem{KO$_1$} Koskela P. and Onninen J.: Mappings on finite distortion: The sharp
modulus of continuity. Trans. Amer. Math. Soc. \textbf{355},
1905--1920 (2003).

\bibitem{KO$_2$} Koskela P., Onninen J.: Mappings of finite distortion: Capacity
and modulus inequalities. J. Reine Agnew. Math. \textbf{599},
1--26 (2006).

\bibitem{KOR} Koskela P., Onninen J. and Rajala K.: Mappings of finite distortion:
Injectivity radius of a local homeomorphism. Future trends in
geometrical function theory, 169--174, Rep. Univ. Jyvaskyla Dep.
Math. Stat., \textbf{92}, Univ. Jyvaskyla, Jyvaskyla (2003).

\bibitem{KR$_*$} Koskela P. and Rajala K.: Mappings of finite distortion: Removable
singularities. Israel J. Math. \textbf{136}, 269--283 (2003).

\bibitem{Koval} Kovalev L.V.: Monotonicity of generalized reduced
modulus. Zapiski Nauch. Sem. POMI \textbf{276}, 219--236 (2001).

\bibitem{KR$_0$} Kovtonyuk D. and Ryazanov V.: On boundaries of space domains. Proc.
Inst. Appl. Math. \& Mech. NAS of Ukraine \textbf{13}, 110--120
(2006) [in Russian].

\bibitem{KR$_1$} Kovtonyuk D. and Ryazanov V.: To the theory of lower
$Q$-homeomorphisms. Ukrainian Math. Bull. \textbf{5} (2), 157--181
(2008).

\bibitem{KR$_2$} Kovtonyuk D. and Ryazanov V.: On the theory of mappings
with finite area distortion. J. d'Anal. Math. \textbf{104},
291--306 (2008).


\bibitem{KR$_3$}
Kovtonyuk D. and Ryazanov V.: On the boundary behavior of
generalized quasi--isometries. ArXiv: 1005.0247 [math.CV], 20 p.
(2010).

\bibitem{KPR} Kovtonyuk D., Petkov I. and Ryazanov V.: On homeomorphisms with finite
distortion in the plane. ArXiv: 1011.3310v2 [math.CV], 1-16 (2010).

\bibitem{Ko} Koronel J.D.: Continuity and $k$-th order differentiability in
Orlicz-Sobolev spaces: $W^kL_A$. Israel J. Math. \textbf{24} (2),
119--138 (1976).

\bibitem{KR} Krasnosel'skii M.A. and Rutitskii Ya.B.: Convex functions and Orlicz spaces.
Noordhoff (1961) (Translated from Russian).

\bibitem{Kr1} Kruglikov V.I.: The existence and uniqueness of mappings that are
quasiconformal in the mean. In: Metric Questions of the Theory of
Functions and Mappings, pp. 123--147. Naukova Dumka, Kiev (1973)
[in Russian].

\bibitem{Kr2} Kruglikov V.I. Capacities of condensors and quasiconformal in the
mean mappings in space. Mat. Sb. \textbf{130} (2), 185--206 (1986)
[in Russian].

\bibitem{KrP} Kruglikov V.I. and Paikov V.I.: Capacities and prime ends of an
$n$-dimensional domain. Dokl. Akad. Nauk Ukrain. SSR Ser. A.
\textbf{5}, 10--13, 84 (1987) [in Russian].

\bibitem{Krush1} Krushkal' S.L.: On mappings quasiconformal in the mean. Dokl.
Akad. Nauk SSSR \textbf{157} (3), 517--519 (1964) [in Russian].

\bibitem{Krush2} Krushkal' S.L.: On the absolute integrability and
differentiability of some classes of mappings of many-dimensional
domains. Sib. Mat. Zh. \textbf{6} (3), 692--696 (1965) [in
Russian].

\bibitem{KK} Krushkal' S.L. and K\"uhnau R.: Quasiconformal mappings, new methods
and applications. Nauka, Novosibirsk (1984) [in Russian].

\bibitem{Kud1} Kud'yavin V.S.: A characteristic property of a class of
$n$-dimensional homeomorphisms. Dokl. Akad. Nauk Ukrain. SSR ser.
A. \textbf{3}, 7--9 (1990) [in Russian].

\bibitem{Kud2} Kud'yavin V.S.: Quasiconformal mappings and $\alpha$-moduli of
families of curves. Dokl. Akad. Nauk Ukrain., \textbf{7}, 11--13
(1992) [in Russian].

\bibitem{Kud3} Kud'yavin V.S.: Estimation of the distortion of distances under
mappings quasiconformal in the mean. Dinam. Sploshn. Sred.,
\textbf{52}, 168--171 (1981) [in Russian].

\bibitem{Kud4} Kud'yavin V.S.: Local boundary properties of mappings
quasiconformal in the mean. In: Collection of Scientific Works,
Institute of Mathematics, Siberian Branch of the Academy of
Sciences of the USSR, pp. 168--171. Novosibirsk (1981) [in
Russian].

\bibitem{Kud5} Kud'yavin V.S.: Behavior of a class of mappings quasiconformal in
the mean at an isolated singular point. Dokl. Akad. Nauk SSSR
\textbf{277} (5), 1056--1058 (1984) [in Russian].

\bibitem{Ku1} K\"uhnau R.: \"Uber Extremalprobleme bei im Mittel quasiconformen
Abbildungen. Lecture Notes in Math. \textbf{1013}, 113--124 (1983)
[in German].

\bibitem{Ku2} K\"uhnau R.: Canonical conformal and quasiconformal mappings.
Identities. Kernel functions. In: Handbook of Complex Analysis,
Geometry Function Theory, vol. 2, pp. 131--163. Elsevier,
Amsterdam (2005).

\bibitem{Ku} Kuratowski K.: Topology, vol. 1. Acad. Press, NY (1968).

\bibitem{Kuz} Kuz'mina G.V.: Moduli of the curve families and quadratic
differentials. Trudy Mat. Inst. AN SSSR \textbf{139}, 1--240 (1980).

\bibitem{LM} Landes R. and Mustonen V.: Pseudo-monotone mappings in Sobolev-Orlicz spaces
and nonlinear boundary value problems on unbounded domains. J.
Math. Anal. Appl. \textbf{88}, 25--36 (1982).

\bibitem{LL} Lappalainen V. and Lehtonen A.: Embedding of Orlicz-Sobolev spaces
in H\"older spaces. Ann Acad. Sci. Fenn. Ser. AI Math. \textbf{14}
(1), 41--46 (1989).

\bibitem{LV} Lehto O. and Virtanen K.: Quasiconformal Mappings in the
Plane. Springer-Verlag, New York (1973).

\bibitem{LF} Lelong-Ferrand J.: Representation conforme et transformations
\`{a} integrale de Dirichlet born\'{e}e. Gauthier-Villars, Paris
(1955).

\bibitem{LSS} Lomako T., Salimov R. and Sevost'yanov E.: On the local behaviour
of homeomorphisms with finite distortion in the plane. ArXiv:
1012.4590v1 [math.CV], 1--18 (2010).

\bibitem{MM} Maly J. and Martio O.: Lusin's condition $(N)$ and mappings of
the class $W^{1,n}$. J. Reine Angew. Math. \textbf{485}, 19--36
(1995).

\bibitem{MV$_1$} Manfredi J.J. and Villamor E.: Mappings with integrable dilatation in
higher dimensions. Bull. Amer. Math. Soc. \textbf{32} (2),
235--240 (1995).

\bibitem{MV$_2$} Manfredi J.J. and Villamor E.: An extension of Reshetnyak's theorem.
Indiana Univ. Math. J. \textbf{47} (3), 1131--1145 (1998).

\bibitem{Ma$_1$} Martio O.: Modern tools in the theory of quasiconformal maps.
Texts in Math. Ser. B, \textbf{27}. Univ. Coimbra, Dept. Mat.,
Coimbra. 1--43 (2000).

\bibitem{MRV} Martio O., Rickman S. and Vaisala J.: Definitions for
quasiregular mappings. Ann. Acad. Sci. Fenn. Ser. A1. Math.
\textbf{448}, 1--40 (1969).

\bibitem{MRSY$_4$} Martio O., Ryazanov V., Srebro U. and Yakubov E.: Mappings with
finite length distortion. J. d'Anal. Math. \textbf{93}, 215--236
(2004).

\bibitem{MRSY$_5$} Martio O., Ryazanov V., Srebro U. and Yakubov E.: $Q$-homeomorphisms.
Contemporary Math. \textbf{364}, 193--203 (2004).

\bibitem{MRSY$_6$} Martio O., Ryazanov V., Srebro U., Yakubov E.: On
$Q$-homeomorphisms. Ann. Acad. Sci. Fenn. \textbf{30}, 49--69
(2005).

\bibitem{MRSY} Martio O., Ryazanov V., Srebro U. and Yakubov E.: Moduli in Modern
Mapping Theory. Springer Monographs in Mathematics. Springer, New
York (2009).

\bibitem{MRV$^*$}
{Martio O., Ryazanov V. and Vuorinen M.}: {BMO and Injectivity of
Space Quasiregular Mappings}. Math. Nachr. \textbf{205}, 149--161
(1999).

\bibitem{MaSa} Martio O. and Sarvas J.: Injectivity theorems in plane and space.
Ann. Acad. Sci. Fenn. Ser. A1 Math. \textbf{4}, 384--401
(1978/1979).

\bibitem{Mattila} Mattila P. Geometry of sets and measures in
Euclidean spaces. Cambridge University Press, Cambridge (1995).

\bibitem{Maz} Maz'ya V.: Sobolev Spaces. Springer-Verlag, Berlin (1985).

\bibitem{Me} Menchoff D.: Sur les differencelles totales des
fonctions univalentes. Math. Ann. \textbf{105}, 75--85 (1931).

\bibitem{Na$_1$} Nakki R.: Boundary behavior of quasiconformal mappings in
$n$-space. Ann. Acad. Sci. Fenn. Ser. A1. Math. \textbf{484},
1--50 (1970).

\bibitem{No} {Norton A.}: {A critical set with nonnull image has large
Hausdorff dimension}. Trans. Amer. Math. Soc. \textbf{296}, no. 1,
367--376 (1986).

\bibitem{On$_1$} Onninen J.: Mappings of finite distortion: future directions and
problems. The $p$-harmonic equation and recent advances in
analysis, Amer. Math. Soc., Contemp. Math., \textbf{370},
199--207, Providence, RI (2005).

\bibitem{On$_2$} Onninen J.: Mappings of finite distortion: Minors of the
differential matrix. Calc. Var. Partial Differ. Eq. \textbf{21}
(4), 335--348 (2004).

\bibitem{On$_3$} Onninen J.: Mappings of finite distortion: Continuity.
Dissertation, Univ. Jyvaskyla, Jyvaskyla, 24 (2002).

\bibitem{Or1} Orlicz  W.: \"Uber eine gewisse Klasse von R\"aumen
vom Typus B. Bull. Intern. de l'Acad. Pol. Serie A, Cracovie
(1932).

\bibitem{Or2} Orlicz W.: \"Uber R\"aume $(L^M)$. Bull. Intern. de l'Acad. Pol. Serie A,
Cracovie (1936).

\bibitem{Pal}
{Palagachev D.K.}: {Quasilinear elliptic equations with VMO
coefficients}. Trans. Amer. Math. Soc. \textbf{347}, no. 7,
2481--2493 (1995).

\bibitem{Pa} Pankka P.: Mappings of finite distortion and weighted
parabolicity. Future trends in geometrical function theory,
175--182, Rep. Univ. Jyvaskyla Dep. Math. Stat., \textbf{92}, Univ.
Jyvaskyla, Jyvaskyla (2003).

\bibitem{Per1} Perovich M.: Isolated singularity of the mean quasiconformal
mappings. Lect. Notes in Math. \textbf{743}, 212--214 (1979).

\bibitem{Per2} Perovich M.: Global homeomorphism of mappings quasiconformal in
the mean. Dokl. Akad. Nauk SSSR \textbf{230} (4), 781--784 (1976).
[in Russian].

\bibitem{Pes} Pesin I.N.: Mappings quasiconformal in the mean. Dokl. Akad. Nauk
SSSR \textbf{187} (4), 740--742 (1969) [in Russian].

\bibitem{Po} Ponomarev S.P.: On the $(N)$-property of homeomorphisms
of the class $W_1^p$. Sibirsk. Mat. Zh. \textbf{28} (2), 140--148
(1987).

\bibitem{QS} {Quinn F. and Sard A.}: {Hausdorff conullity of critical
images of Fredholm maps}. Amer. J. Math. \textbf{94}, 1101--1110
(1972).

\bibitem{RR$^*$} Rado T., and Reichelderfer P.V.: Continuous Transformations in Analysis.
Springer-Verlag, Berlin (1955).

\bibitem{Ra} {Ragusa M.A.}:
{Elliptic boundary value problem in vanishing mean oscillation
hypothesis}. Comment. Math. Univ. Carolin. \textbf{40}, no. 4,
651--663 (1999).

\bibitem{Ra$_1$} Rajala K.: Mappings of finite distortion: the Rickman--Picard
theorem for mappings of finite lower order. J. Anal. Math.
\textbf{94}, 235--248 (2004).

\bibitem{Ra$_2$} Rajala K.: Mappings of finite distortion: removability of Cantor
sets. Ann. Acad. Sci. Fenn. Math. \textbf{29} (2), 269--281 (2004).

\bibitem{Ra$_3$} Rajala K.: Mappings of finite distortion: removable singularities
for locally homeomorphic mappings. Proc. Amer. Math. Soc.
\textbf{132} (11), 3251--3258 (2004) (electronic).

\bibitem{Ra$_4$} Rajala K.: Mappings of finite distortion: Removable singularities.
Dissertation, Univ. Jyvaskyla, Jyvaskyla, pp. 74 (2003).

\bibitem{RZZ} Rajala K., Zapadinskaya A. and Z\"urcher T.:
Generalized Hausdorff dimension distortion in euclidean spaces under
Sobolev mappings. ArXiv:1007.2091v1 [math.CA], 1--13 (2010).

\bibitem{Rei} Reimann H.M.: On the absolute continuity of surface representation.
Comment. Math. Helv. \textbf{46}, 44--47 (1971).

\bibitem{ReRy} Reimann H.M. and Rychener T.: Funktionen Beschr\"ankter Mittlerer
Oscillation. Lecture Notes in Math. \textbf{487} (1975).

\bibitem{Re} Reshetnyak Yu.G.: Space mappings with bounded distortion. Nauka,
Novosibirsk (1982); English transl., Translations of Mathematical
Monographs, vol. 73, Amer. Math. Soc., Providence, RI (1988).

\bibitem{Resh} Reshetnyak Yu.G.: The condition $(N)$ for $W^{1,n}_{\rm loc}$
space mappings. Sibirsk. Math. Zh. \textbf{28}, 149--153 (1987) [in
Russian].

\bibitem{Re$_1$} Reshetnyak Yu.G.: Sobolev classes of functions with
values in a metric space. Sibirsk. Mat. Zh. \textbf{38}, 657--675
(1997).

\bibitem{Re$_2$} Reshetnyak Yu.G.: Some geometric properties of functions and
mappings with generalized derivatives. Sibirsk. Mat. Zh. \textbf{7},
886--919 (1966).

\bibitem{Ri} Rickman S.: Quasiregular Mappings. Springer, Berlin etc.
(1993).

\bibitem{Rya1} Ryazanov V.I.: On compactification of classes with integral
restrictions on the Lavrent'ev characteristics. Sibirsk. Mat. Zh.
\textbf{33} (1), 87--104 (1992) [in Russian]; translation in
Siberian Math. J. \textbf{33} (1), 70--86 (1992).

\bibitem{Rya2} Ryazanov V.I.: On quasiconformal mappings with measure
restrictions. Ukrain. Mat. Zh. \textbf{45} (7), 1009--1019 (1993)
[in Russian]; translation in Ukrain. Math. J. \textbf{45} (7),
1121--1133 (1993).

\bibitem{Rya3} Ryazanov V.I.: On mappings that are quasiconformal in the mean.
Sibirsk. Mat. Zh. \textbf{37}(2), 378--388 (1996) [in Russian];
translation in Siberian Math. J. \textbf{37} (2), 325--334 (1996).

\bibitem{RS} Ryazanov V. and Sevost'yanov E.: Toward the theory of ring
$Q$-homeomorphisms. Israel J. Math. \textbf{168}, 101--118 (2008).

\bibitem{RS_*} Ryazanov V. and Sevost'yanov E.: Equicontinuity of
mappings quasiconformal in the mean. ArXiv: 1003.1199v4 [math.CV],
1-16 (2010).

\bibitem{RSY$_1$} Ryazanov V., Srebro U. and Yakubov E.: On ring
solutions of Beltrami equation. J. d'Anal. Math. \textbf{96},
117--150 (2005).

\bibitem{RSY$_2$} Ryazanov V., Srebro U. and Yakubov E.:
Integral conditions in the theory of the Beltrami equations. Complex
Variables and Elliptic Equations (to appear).

\bibitem{Sa} Saks S.: Theory of the Integral. Dover, New York (1964).

\bibitem{Sar} Sard A.: The measure of the critical values of
differentiable maps. Bull. Amer. Math. Soc. \textbf{48}, 883--890
(1942).

\bibitem{Sard} Sard A.: The equivalence of $n$-measure and Lebesgue measure in $E_n$.
Bull. Amer. Math. Soc. \textbf{49}, 758--759 (1943).


\bibitem{Sard$_2$} {Sard A.}: {Images of critical sets}. Ann. Math.  \textbf{68}, no. 2,
247--259 (1958).

\bibitem{Sard$_3$} {Sard A.}: {Hausdorff measure of
critical images on Banach manifolds}. Amer. J. math. \textbf{87},
158--174 (1965).

\bibitem{Shl} Shlyk V.A.: The equality between $p$-capacity and
$p$-modulus. Sib. Math. J. \textbf{34} (6), 1196--1200 (1993).

\bibitem{So} Sobolev S.L.: Applications of functional analysis in mathematical
physics. Izdat. Gos. Univ., Leningrad (1950); English transl. Amer.
Math. Soc., Providence, R.I. (1963).

\bibitem{Sol} Solynin A.Yu.: Moduli and extremal-metric problems.
Algebra and Analysis \textbf{11} (1), 3--86 (1999).

\bibitem{Str} Strugov Y.F.: Compactness of the classes of mappings
quasiconformal in the mean. Dokl. Akad. Nauk SSSR \textbf{243}
(4), 859--861 (1978) [in Russian].

\bibitem{StrSy} Strugov Y.F. and Sychov A.V.: On different classes of mappings
quasiconformal in the mean. Vest. PANI, \textbf{7}, 14--19 (2002)
[in Russian].

\bibitem{Su1} Suvorov G.D.: Generalized principle of the length and area in the
mapping theory. Naukova Dumka, Kiev (1985) [in Russian].

\bibitem{Su2} Suvorov G.D.: The metric theory of prime ends and boundary
properties of plane mappings with bounded Dirichlet integrals.
Naukova Dumka, Kiev (1981) [in Russian].

\bibitem{Su3} Suvorov G.D.: Families of plane topological mappings.
AN SSSR, Novosibirsk (1965) [in Russian].

\bibitem{Sych} Sychev A.V.: Moduli and spatial quasiconformal
mappings. Novosibirsk, Nauka (1983).

\bibitem{Tu} Tuominen H.: Characterization of Orlicz-Sobolev space.
Ark. Mat. \textbf{45} (1), 123--139 (2007).

\bibitem{UV} Ukhlov A. and Vodop'yanov S.:
{Sobolev spaces and $(P,Q)-$quasiconformal mappings of Carnot
groups}. Siberian Math. J. \textbf{39}, 665--682 (1998).

\bibitem{UkhVo} Ukhlov A. and Vodop'yanov S.: Mappings associated with weighted
Sobolev Spaces. Complex Anal. Dynam. Sys. III. Contemp. Math.
\textbf{455}, 363--382 (2008).

\bibitem{Va$_1$} V\"ais\"al\"a J.: Lectures on $n$-dimensional quasiconformal
mappings. Lecture Notes in Math. \textbf{229}, Springer-Verlag,
Berlin (1971).

\bibitem{Va$_2$} V\"ais\"al\"a J.: On quasiconformal mappings in space.
Ann. Acad. Sci. Fenn. Ser. A1. Math. \textbf{298}, 1--36 (1961).

\bibitem{Va$_3$} V\"ais\"al\"a J.: On the null-sets for extremal distances. Ann.
Acad. Sci. Fenn. Ser. A1. Math. \textbf{322}, 1--12 (1962).

\bibitem{Vas} Vasil'ev A.: Moduli of families of curves for conformal and
quasiconformal mappings. Lecture Notes in Math. (1788),
Springer-Verlag, Berlin–New York (2002).

\bibitem{Vo} Vodop'yanov S.: Mappings with bounded distortion and with finite distortion on
Carnot groups. Sibirsk. Mat. Zh. \textbf{40} (4), 764--804 (1999);
transl. in Siberian Math. J. \textbf{40} (4), 644--677 (1999).

\bibitem{Vu} Vuillermot P.A.: H\"older-regularity for the solutions of strongly nonlinear
eigenvalue problems on Orlicz-Sobolev spaces. Houston J. Math.
\textbf{13}, 281--287 (1987).

\bibitem{Wi} Wilder R.L.: Topology of Manifolds. AMS, New York (1949).

\bibitem{Wh} Whitney H.: A function not constant on a connected set
of critical points. Duke Math. J. \textbf{1}, 514--517 (1935).

\bibitem{Za} Zaanen A.C.: Linear analysis. Noordhoff (1953).

\bibitem{Zi} Ziemer W.P.: Extremal length and conformal capacity. Trans. Amer.
Math. Soc. \textbf{126} (3), 460--473 (1967).

\bibitem{Zo1} Zorich V.A.: Admissible order of growth of the characteristic of
quasiconformality in the Lavrent'ev theory. Dokl. Akad. Nauk SSSR
\textbf{181} (1968) [in Russian].

\bibitem{Zo2} Zorich V.A.: Isolated singularities of mappings with bounded
distortion. Mat. Sb. \textbf{81}, 634--638 (1970) [in Russian].

\bibitem{Zo3} Zorich V.A.: Quasiconformal mappings and the asymptotic
geometry of manifolds. Russ. Math. Surv. \textbf{57} (3), 437--462
(2002); transl. from Usp. Mat. Nauk \textbf{57} (3), 3--28 (2002).

\end{thebibliography}
\end{document}